\documentclass[12pt]{article}
\usepackage[utf8]{inputenc}
\usepackage[margin=1in]{geometry}
\usepackage{amsmath,amsthm,amssymb,amsfonts, pifont, bbm, scrextend,}
\usepackage[usenames,dvipsnames]{xcolor}
\usepackage{fancyhdr}
\usepackage{array}
\usepackage {tikz}
\usetikzlibrary {positioning}
\definecolor {processblue}{cmyk}{0.96,0,0,0}
\definecolor {processred}{cmyk}{0,1,.5,.25}
\definecolor{forestgreen}{RGB}{0,98,51}
\usepackage{appendix}
\usepackage{fancyhdr}
\usepackage{color,colortbl}
\newcounter{CommentCounter}

\definecolor{Gray}{gray}{0.9}
\definecolor{LightGreen}{rgb}{0,0.1,0.1}
\newcommand{\upto}{\!\nearrow\!}

\newcommand{\R}{\mathbb R}
\newcommand{\Gr}{\mathcal G}
\newcommand{\Hg}{\mathcal H}
\newcommand{\Pa}{{P}}
\newcommand{\Pm}{{\bf P}}
\newcommand{\Em}{{\bf E}}
\newcommand{\Ea}{{E}}

\newcommand{\N}{\mathbb N}

\newcommand{\Z}{\mathbb Z}
\newcommand{\0}{\vec 0}

\newcommand{\given}{\hspace{3pt}\vline\hspace{3pt}}

\newcommand{\cpath}{\overline{\pi}}
\newcommand{\vecuk}{u}
\newcommand{\vecu}{u}
\newcommand{\vecvk}{v}
\newcommand{\vecv}{v}

\usepackage{tikz}
\usepackage{pgfplots}
\usepackage{eqnarray}
\usepackage{outlines}
\usepackage[mathscr]{euscript}
 \let\mathscr\relax
\usepackage[scr]{rsfso}
\usepackage{multicol}
\usepackage{enumerate}
\usepackage[colorlinks=true, pdfstartview=FitV, linkcolor=blue,
citecolor=blue, urlcolor=blue]{hyperref}

\newtheorem{thm}{Theorem}[section]
\newtheorem*{thm*}{Theorem}
\numberwithin{thm}{section}

\newtheorem{prop}[thm]{Proposition}
\newtheorem{lem}[thm]{Lemma}
\newtheorem{conj}[thm]{Conjecture}

\numberwithin{quest}{section}

\theoremstyle{remark}

\numberwithin{rem}{section}

\theoremstyle{definition}

\newtheorem{claim}{Claim}
\numberwithin{claim}{exmp}

\pgfplotsset{compat=1.12}

\begin{document}
\author{
Daniel J. Slonim\footnote{Purdue University,
  Department of Mathematics,
150 N. University Street, West Lafayette, IN 47907, dslonim@purdue.edu}
}
\title{Directional transience of random walks in random environments with bounded jumps}

\date{\today}

\maketitle

\begin{abstract}
    This paper has two main results, which are connected through the fact that the first is a key ingredient in the second. Both are extensions of results concerning directional transience of nearest-neighbor random walks in random environments to allow for bounded jumps.
    Zerner and Merkl \cite{Zerner&Merkl2001} proved a 0-1 law for directional transience for planar random walks in random environments. We extend the result to non-planar i.i.d. random walks in random environments on $\Z^2$ with bounded jumps. 
    Sabot and Tournier \cite{Sabot&Tournier2016} characterized directional transience for a given direction for nearest-neighbor random walks in Dirichlet environments on $\Z^d$, $d\geq1$. We extend this characterization to random walks in Dirichlet environments with bounded jumps.
    
    \\
    {\it Kewords:}
    random walk,
    random environment,
    bounded jumps,
    0-1 law,
    Dirichlet environments,
    directional transience,
    recurrence
\end{abstract}
\section{Introduction}

In 1981, Kalikow \cite{Kalikow1981} asked whether, for i.i.d. random walks in random environments (RWRE) in 2 dimensions, the $x$-coordinate of the walker's position must approach infinity with probability either 0 or 1. In 2001, Zerner and Merkl \cite{Zerner&Merkl2001} answered this question in the affirmative for nearest-neighbor, i.i.d., elliptic RWRE in 2 dimensions, and not just for the horizontal component but for the component in any direction $\ell\in S^1$, where $S^1$ is the unit circle in $\R^2$. (They also showed that the i.i.d. assumption is necessary). Such a 0-1 law is still an open conjecture for dimensions $d\geq3$.
In this paper, we extend the result of Zerner and Merkl by removing the nearest-neighbor assumption, showing that for i.i.d. elliptic RWRE with bounded jumps on $\Z^2$, the 0-1 law holds for all directions $\ell\in S^1$. Our approach is largely based on that of \cite{Zerner2007}, which is a simplification of the proof given in \cite{Zerner&Merkl2001}. However, the removal of the nearest-neighbor assumption creates a need for additional work and a number of adjustments.

We then turn our attention to random walks in Dirichlet environments (RWDE). For a given direction $\ell\in S^{d-1}$ (where $S^{d-1}$ is the unit sphere in $\R^d$), the question of transience and recurrence in direction $\ell$ is completely understood for nearest-neighbor RWDE. Tournier remarks in \cite{Tournier2015} that many of the results used in the characterization of transience do not rely on the nearest-neighbor assumption, so much of what was known in the nearest-neighbor case carries over to the bounded-jumps case. However, not everything carries over directly. One crucial step toward characterizing directional transience is a 0-1 law. As Tournier points out in his remark, the proof of the 0-1 law for RWDE in dimensions $d\geq3$ given in \cite{Bouchet2013} does not require the nearest-neighbor assumption, but the proof for RWRE in dimension $d=2$ in \cite{Zerner&Merkl2001} and \cite{Zerner2007} does require the nearest-neighbor assumption. Our extension of the 0-1 law for $d=2$ to bounded jumps means that for RWDE with bounded jumps, the 0-1 law is now proven for all dimensions. Removing the nearest-neighbor assumption creates one other obstacle to fully characterizing directional transience in a given direction. When the annealed drift is zero, the nearest-neighbor argument relies on a symmetry that does not necessarily exist in the bounded-jump case. Therefore, completing a characterization of directional transience amounts to showing, for all dimensions, that zero annealed drift implies recurrence in any direction. 

The rest of this section formally defines our model. Section \ref{sec:main} states and proves the 0-1 law. Section \ref{sec:Dirichlet} copmletes the characterization of directional transience for a given direction in the Dirichlet case.

\subsection{Model}
Let $V$ be a finite or countable set. 
For the theorems stated in this paper, we always have $V=\Z^d$ (with $d=2$ in Section \ref{sec:main}), but the proofs in Section \ref{sec:Dirichlet} will require constructing RWRE on other sets. 
An {\em environment on $V$} is a nonnegative function $\omega:V\times V\to [0,1]$ such that for all $x\in V$, $\sum_{y\in V}\omega(x,y)=1$. We denote by $\Omega_V$ the set of all environments on $V$. For a given environment $\omega$ and $x\in V$, we can define the {\em quenched} measure $P_{\omega}^x$ on $V^{\N}$ (where we assume $0\in \N$) to be the law of a Markov chain ${\bf X}=(X_n)_{n\geq0}$ on $V$, started at $x$, with transition probabilities given by $\omega$. That is: $P_{\omega}^x(X_0=x)=1$, and for $n\geq1$, $P_{\omega}^x(X_{n+1}=y|X_0,\ldots,X_n)=\omega(X_n,y)$. 

Let $\mathcal{F}_V$ be the Borel sigma field with respect to the product topology on $\Omega_V$,
and let $\Pm$ be a probability measure on $(\Omega_V,\mathcal{F}_V)$. For a given $x\in V$, we define the {\em annealed} measure $\Pa^x=\Pm\times P_{\omega}^x$ on $\Omega_V\times V^{\N}$ by
\begin{equation*}\Pa^x(A\times B)=\int_{A}P_{\omega}^x(B)\Pm(d\omega)\end{equation*}
for measurable $A\subset\Omega_V,B\subset V^{\N}$. In particular, for measurable $B\subset V^{\N}$, $\Pa^x(\Omega_V\times B)=\Em[P_{\omega}^x(B)]$. We often abuse notation by writing $\Pa^x(B)$ instead of $\Pa^x(\Omega_V\times B)$. 
When referring to ``the law'' of a RWRE, we mean the annealed law unless otherwise specified.

All measures $\Pm$ on $\Omega_{\Z^d}$ considered in this paper satisfy the following conditions.
\begin{enumerate}[(C1)]
    \item\hypertarget{cond:C1} Under $\Pm$, if $\omega^x(y)=\omega(x,x+y)$, the 
    $(\omega^x)_{x\in\Z^d}$ are i.i.d.;
    \item\hypertarget{cond:C2} With $\Pm$-probability 1, the Markov chain induced by $\omega$ has only one infinite communicating class, and it is reachable from every site.
    \item\hypertarget{cond:C3} There is an $R>0$ such that with $\Pm$-probability 1, $\omega(x,y)=0$ whenever $|x-y|>R$;
\end{enumerate}

Condition (\hyperlink{cond:C2}{C2}) replaces the weak ellipticity assumption \cite{Zerner&Merkl2001} and \cite{Zerner2007}, under which $\Pm(\omega(0,y)>0)=1$ for all four nearest neighbors $y$ of 0. Our condition (\hyperlink{cond:C2}{C2}) is satisfied whenever the Markov chain induced by $\omega$ is $\Pm$--almost surely irreducible. In particular, if there is a set of possible jumps that always have positive probability for $\Pm$--almost every environment, and it is possible to reach any site from any other site using such jumps, then  (\hyperlink{cond:C2}{C2}) is satisfied (under the weak ellipticity assumption of \cite{Zerner&Merkl2001} and \cite{Zerner2007}, the set of nearest-neighbor jumps has this property). For an example where the Markov chain is $\Pm$--almost surely \textit{not} irreducible, but condition (\hyperlink{cond:C2}{C2}) is still satisfied, see Appendix \ref{append:C2}.

For $\ell\in S^{d-1}$, define
$$
A_{\ell}:=\{{\bf X}\in(\Z^d)^{\N}:\lim_{n\to\infty}X_n\cdot\ell=\infty\}.
$$

\section{The 0-1 Law for Dimension 2}\label{sec:main}

We prove the 0-1 law for directional transience for i.i.d. RWRE on $\Z^2$ with bounded jumps.

\begin{thm}\label{thm:01Law}
Let $d=2$, let assumptions (\hyperlink{cond:C1}{C1}), (\hyperlink{cond:C2}{C2}), and (\hyperlink{cond:C3}{C3}) hold, and let $\ell\in S^{1}$. Then $\Pa^0(A_{\ell})\in\{0,1\}$.
\end{thm}

Before giving the proof, we summarize  Zerner's proof in \cite{Zerner2007} and discuss where ours will differ. The idea of the proof from \cite{Zerner2007} is that if the probability of transience in both direction to the left and to the right (for instance) is positive, then with non-vanishing probability, one should be able to start two walks in the same environment on different sides of a wide strip and have both walks cross the strip and exit on the opposite side from where they started (call this the {\em strip traversal event}). If the starting points are chosen correctly, this should lead to the paths of the walks crossing at least half the time. The nearest-neighbor assumption is leveraged here, as it implies that crossing paths must intersect.
This intersection entails a low-probility event. Consider the meeting point. A walk came a long way from the right to hit that point, so it must have a very high probability of transience to the left. But another walk went through this same point, and then traveled a long distance to the right. The probability of such an event can be made arbitrarily low, implying that transience to the left or to the right must have zero probability. 

We now discuss differences between the above argument and ours. The argument from \cite{Zerner2007} breaks the strip-traversal event into the event that the two paths intersect and the event that they do not intersect, the latter being a subset of the event that the walks land on opposite sides of a straight line through their starting points. We must consider three events: that the paths intersect, that the paths come within a specified distance of each other, and that the paths land on opposite sides of the line through their starting points. Showing that these events are the entirety of the strip-traversal event is a step not required in Zerner's argument. If the linear interpolations of the paths cross, then the finite range assumption implies that the paths must come near each other. However, unlike in the planr case, it is possible for them to land on the same side of the line through their starting points without the linear interpolations crossing, and we must show that this event also entails the walks coming near each other.

A more significant difficulty is in comparing the probability of the event that the walks come near each other to the probability that they actually meet. Because we are not assuming uniform ellipticity, and because a part of the environment where both paths come near each other is not 
necessarily a ``typical'' part of the environment, arguments dealing with the quenched probability of a modified path that causes an intersection would be difficult. Instead, we focus on annealed probabilities, which requires us to leverage independence of the environment at different sites while still forcing the walks to meet. This requires careful attention to the work of defining the right stopping times and events, and unlike in Zerner's argument, results in our defining a meeting event where one of the walks does not necessarily complete the strip traversal, but which nonetheless has vanishing probability.

\begin{proof}[Proof of Theorem \ref{thm:01Law}]
We divide this proof into steps.

\hspace{-.25in}{\em Step 1: Preliminaries}

Fix $\ell\in S^1$. 
For $a\in\R$
and ${\bf X}$ a 
sequence in $\Z^2$, consider the stopping times
$$
T_{\geq a}=T_{\geq a}({\bf X}):=\inf\{n\geq0:(X_n\cdot\ell)\geq a\},
$$
and likewise for $T_{\leq a}$, $T_{>a}$, and $T_{<a}$.
Similarly, for a set $S\subset{\Z^2}$, define
$$
T_S=T_{S}({\bf X}):=\inf\{n\geq0:X_n\in S\}.
$$
We often suppress the argument ${\bf X}$ when the sequence intended is clear from the context.

For $y\in \Z^2$ and a path $\gamma=(x_0,x_1,\ldots,x_n)$ (which is a path of length $n$, and is a loop if $x_n=x_0$), define $y+\gamma:=(y+x_0,y+x_1,\ldots,y+x_n)$. This is simply a space shift of the path $\gamma$. 
For a path $\gamma=(x_0,x_1,\ldots,x_n)$, we will talk about the {\em annealed probability of $\gamma$}, or the {\em probability that ${\bf X}$ takes $\gamma$}. This simply means 
$$
\Pa^{x_0}(X_0=x_0,X_1=x_1,\ldots,X_n=x_n).
$$
Call a path $\gamma$ a {\em possible path} if it has positive annealed probability. Note that a loop $(x_0)$ of zero length has annealed probability 1, since $\Pa^{x_0}(X_0=x_0)=1$.

Now by assumption (\hyperlink{cond:C2}{C2}), there is a possible path connecting any two points. This is because if $x,y\in\Z^d$, then with positive $\Pm$-probability, $y$ is in the infinite communicating class and thus reachable from $x$ by some finite path. Since there are countably many finite paths from $x$ to $y$, at least one must have positive annealed probability. Let $M$ be large enough that for any vertex $y$ in a closed unit disc of radius $2R$ centered at 0, there is a path of positive probability from $0$ to $y$ with length no more than $M$.

By Kalikow's 0-1 law (Theorem \ref{thm:K01Law} in Appendix \ref{app:Kalikow}), $\Pa^0(A_{\ell}\cup A_{-\ell})\in\{0,1\}$. Thus, it suffices to show that $\Pa^0(A_{\ell})\Pa^0(A_{-\ell})=0$ under the assumption that $\Pa^0(A_{\ell}\cup A_{-\ell})=1$. 
Lemma \ref{lem:986} tells us that
$\Pa^0(A_{\ell})>0$ if and only if $\Pa^0(T_{<0}=\infty)>0$ and $\Pa^0(A_{-\ell})>0$ if and only if $\Pa^0(T_{>0}=\infty)>0$. Thus, it suffices to show 
\begin{equation}\label{eqn:222}
\Pa^0(T_{<0}=\infty)\Pa^0(T_{>0}=\infty)=0.
\end{equation}

For $a,b\in\R$, define the event
\begin{equation}\notag
    G_a^b:=\begin{cases}
    \vspace{-.19in}
    \\
    \left\{T_{\geq b}<T_{<a}\right\}\quad\quad\text{if }b>a;
    \\
    \vspace{-.15in}
    \\
    \left\{T_{\leq b}<T_{>a}\right\}\quad\quad\text{if }b<a
    \\
    \vspace{-.19in}
    \end{cases}.
\end{equation}
Note that for fixed $a$, 
\begin{equation}\label{eqn:243}
\lim_{b\to\infty}G_a^b=\bigcap_{b>a}G_a^b\subset\{T_{<a}=\infty\};
\quad
\lim_{b\to-\infty}G_a^b=\bigcap_{b<a}G_a^b\subset\{T_{>a}=\infty\}. 
\end{equation}

\hspace{-.25in}{\em Step 2: Two walks in one environment}

In this step, we define a point $z_L$ and measures on $(\Z^2)^{\mathbb{N}}\times(\Z^2)^{\N}$ that encapsulate the notion of running two random walks in the same environment, one from $0$ and one from $z_L$. We define a ``strip-traversal event'' in which the two walks cross a strip in opposite directions with certain restrictions, and then define three subsets of this event and show that their union is the whole event. Showing this last statement requires more delicate work than showing the analagous statement in \cite{Zerner2007}, due to fringe cases that are depicted in figure \ref{fig:PossibleSituations}, which do not appear in the nearest-neighbor model.

Fix a unit vector $\ell^{\perp}$ perpendicular to $\ell$. Choose a sequence $z_L\in\Z^2$ indexed by $L\in\N$ such that 
\begin{itemize}
    \item $z_L\cdot\ell\geq2L$,
    \item With positive $\Pa^0$-probability, $X_{T_{\geq2L}}=z_L$, and
    \item $z_L\cdot\ell^{\perp}$ is a median of the distribution of $X_{T_{\geq2L}}\cdot\ell^{\perp}$ under the measure $\Pa^0(\cdot|G_0^{2L})$. That is, $\Pa^0(X_{T_{\geq2L}}\cdot\ell^{\perp} >z_L\cdot\ell^{\perp}|G_0^{2L})\leq\frac12$ and $\Pa^0(X_{T_{\geq2L}}\cdot\ell^{\perp} <z_L\cdot\ell^{\perp}|G_0^{2L})\leq\frac12$.
\end{itemize}
Define $x_L:=z_L\cdot\ell$.
Due to the allowance of jumps, $z_L$ may not be uniquely defined for each $L$---for example, if  $\ell=(1,0)$ and a jump of two steps to the right is possible, then $(2L,h)$ and $(2L+1,h)$ would both be candidates for $z_L$ for some $h$---but one may, for instance, always take the candidate with the smallest $\ell$ component.  Now consider two independent random walks ${\bf X}^1=(X_n^1)_n$ and ${\bf X}^2=(X_n^2)_n$ moving in the same environment, with the first walk starting at $0$ and the second starting at $z_L$. For $\omega\in\Omega_{\Z^2}$ and $a,b\in\Z^2$, let $P_{\omega}^{a,b}$ be the product measure $P_{\omega}^a\times P_{\omega}^b$ on the set $(\Z^2)^{\N}\times(\Z^2)^{\N}$ with typical element $({\bf X}^1,{\bf X}^2)$. Let $\Pa^{a,b}$ be the corresponding annealed measure. 

We consider the ``strip traversal event'' $G_0^{2L}\times G_{x_L}^0$, which is roughly the event that both walks cross the strip $\{0\leq x\cdot\ell\leq2L\}$ before leaving it; the walk starting at 0 is in $G_0^{2L}$, while the walk starting at $z_L$ is in $G_{x_L}^0$. Zerner shows\footnote{Zerner actually shows $\Pa^0(T_{<0}=\infty)\Pa^0(T_{>0}=\infty)\leq\liminf_{L\to\infty}\Pa^{0,z_L}(G_0^{2L}\times G_{x_L}^0)$, and this is all that is needed for his argument and ours. However, a very brief and straightforward addition to Zerner's argument would show that the liminf is actually a limit and that equality holds, so we write it that way for cosmetic reasons.
} \cite[equation (10)]{Zerner2007} that
\begin{equation}\label{eqn:395}
    \Pa^0(T_{<0}=\infty)\Pa^0(T_{>0}=\infty)=\lim_{L\to\infty}\Pa^{0,z_L}(G_0^{2L}\times G_{x_L}^0).
\end{equation}

Now consider the following three subsets of the strip traversal event: 

\begin{itemize}
    \item $\mathcal{O}_L$, the {\em opposite-sides event}. This is the event that ${\bf X}^1\in G_0^{2L}$, ${\bf X}^2\in G_{x_L}^0$, and $\left[(X_{T_{\geq2L}}^1-z_L)\cdot\ell^{\perp}\right]\left[X_{T_{\leq0}}^2\cdot\ell^{\perp}\right]<0$. 
    \item $\mathcal{I}_L$, the {\em intersection event}. This is the event that ${\bf X}^1\in G_0^{2L}$, ${\bf X}^2\in G_{x_L}^0$,
    and for some $0\leq m \leq T_{\geq 2L}({\bf X}^1)$, $0\leq n \leq T_{\leq0}({\bf X}^2)$, $X_m^1=X_n^2$.
    \item $\mathcal{P}_L$, the {\em proximity event}. This is the event that ${\bf X}^1\in G_0^{2L}$, ${\bf X}^2\in G_{x_L}^0$, 
    and for some $0\leq m \leq T_{\geq 2L}({\bf X}^1)$, $0\leq n \leq T_{\leq0}({\bf X}^2)$, $|X_m^1-X_n^2|\leq2R$.
\end{itemize}
Clearly $\mathcal{I}_L\subset\mathcal{P}_L$. We claim that the three events together comprise the entirety of the strip traversal event.

\setcounter{exmp}{1}
\begin{claim}\label{claim:WholeEvent}
\begin{equation}\label{eqn:282}
G_0^{2L}\times G_{x_L}^0
=
\mathcal{O}_L\cup\mathcal{P}_L
=
\mathcal{O}_L\cup\mathcal{P}_L\cup\mathcal{I}_L.
\end{equation}
\end{claim}

The events $\mathcal{O}_L$ and $\mathcal{P}_L$ are each specified to be contained in the event $G_0^{2L}\times G_{x_L}^0$, so their union is as well. Now assume ${\bf X}^1\in G_0^{2L}$ and ${\bf X}^2\in G_{x_L}^0$. We will show that either $\mathcal{O}_L$ or $\mathcal{P}_L$ occurs.
Let $\cpath^1$ be the continuous linear interpolation of the path taken by ${\bf X}^1$, and let $\cpath^2$ be the continuous linear interpolation of the path taken by ${\bf X}^2$. Let $\alpha_2$ be the last point in $\R^2$ where $\cpath^2$ crosses the line $\{x\cdot\ell=2L\}$. Let $\beta_1$ be the first point where $\cpath^1$ crosses $\{x\cdot\ell=2L\}$, and let $\beta_2$ be the first point where $\cpath^2$ crosses $\{x\cdot\ell=0\}$. Let $z_L'$ be the point on the line $\{x\cdot\ell=2L\}$ with $(z_L-z_L')\cdot\ell^{\perp}=0$ (thus, $z_L=z_L'+(x_L-2L)\ell$). Note $\alpha_2$, $\beta_1$, $\beta_2$, and $z_L'$ need not be in $\Z^d$.

To show that either $\mathcal{O}_L$ or $\mathcal{P}_L$ occurs, we will assume $\mathcal{O}_L$ does not occur and prove that $\mathcal{P}_L$ must occur. If $\mathcal{O}_L$ does not occur, then $(X_{T_{\geq2L}}^1-z_L)\cdot\ell^{\perp}$ and $X_{T_{\leq0}}^2\cdot\ell^{\perp}$ are either both positive or both negative, or else at least one is 0. If $(X_{T_{\geq2L}}^1-z_L)\cdot\ell^{\perp}=0$, then $\mathcal{P}_L$ occurs, because $X_{T_{\geq2L}}$ and $z_L$ have the same $\ell^{\perp}$ component and both have $\ell$ component between $2L$ and $2L+R$. Similarly, if $X_{T_{\leq0}}^2\cdot\ell^{\perp}=0$, then $\mathcal{P}_L$ occurs. 

Now suppose $(X_{T_{\geq2L}}^1-z_L)\cdot\ell^{\perp}$ and $X_{T_{\leq0}}^2\cdot\ell^{\perp}$ are both nonzero and have the same sign. Without loss of generality, we may assume both are positive (otherwise, rename the directions $\ell^{\perp}$ and $-\ell^{\perp}$). To show that $\mathcal{P}_L$ occurs, we must show that for some $0\leq m\leq T_{\geq 2L}({\bf X}^1)$ and for some $0\leq n < T_0^2$, $|X_m^1-X_n^2|\leq2R$.

First, suppose that $\beta_2\cdot\ell^{\perp}<0$; this situation is depicted in on the left in figure \ref{fig:PossibleSituations}. Then $X_{T_{\leq0}-1}^2\cdot\ell>0$ and $X_{T_{\leq0}-1}^2\cdot\ell^{\perp}<0$, but $X_{T_{\leq0}}^2\cdot\ell<0$ and $X_{T_{\leq0}}^2\cdot\ell^{\perp}>0$. In one step, the walker that started at $z_L$ crosses the line $\{x\cdot\ell=0\}$ and the line $\{x\cdot\ell^{\perp}=0\}$. It follows that $X_{T_{\leq0}-1}^2$ must be within a radius $R$ of 0, and the event $\mathcal{P}_L$ occurs. 
\begin{figure}
    \centering
\includegraphics[height=1.5in]{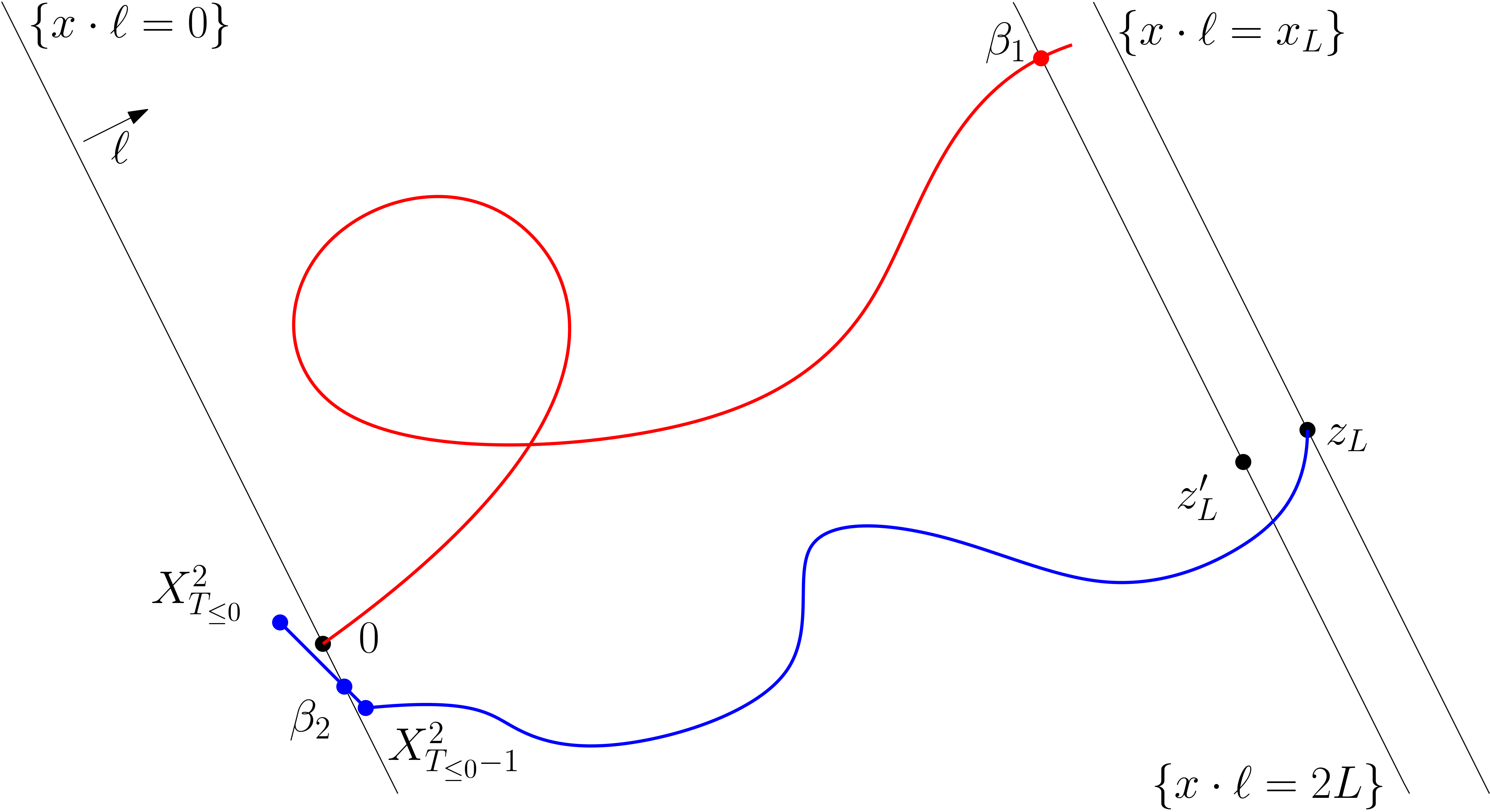}
\hspace{.3in}
\includegraphics[height=1.5in]{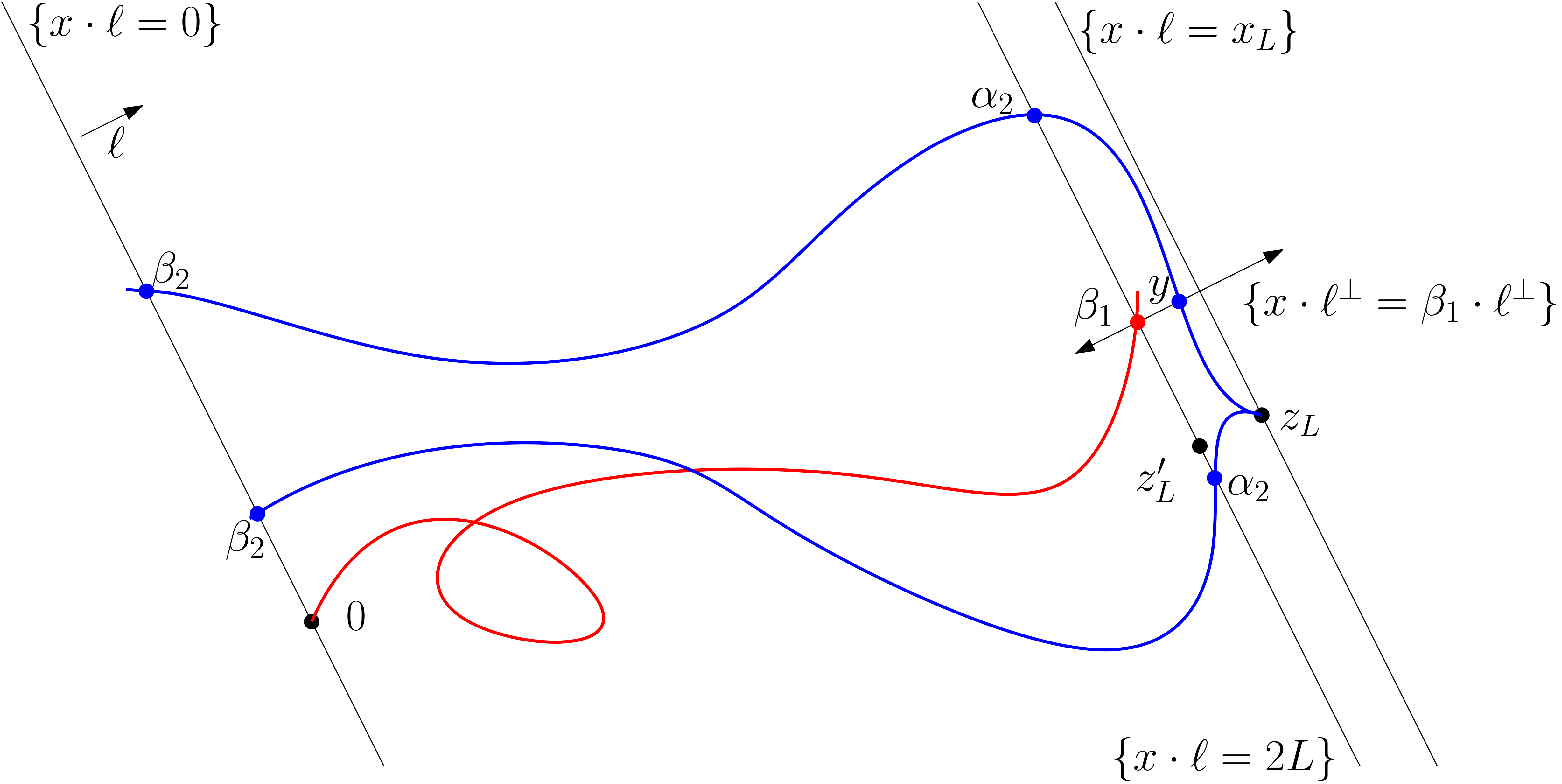}
    \caption{On the left, we have $(X_{T_{\geq2L}}^1-z_L)\cdot\ell^{\perp}>0$ and $X_{T_{\leq0}}^2\cdot\ell^{\perp}>0$, but $\beta_2\cdot\ell^{\perp}<0$.
    On the right,  different ways to have $\beta_2\cdot\ell^{\perp},~ (\beta_1-z_L')\cdot\ell^{\perp}>0$ are depicted. The upper path from $z_L$ shows the situation $\alpha_2\cdot\ell^{\perp}>\beta_1\cdot\ell^{\perp}$, while the lower path from $z_L$ shows $\alpha_2\cdot\ell^{\perp}<\beta_1\cdot\ell^{\perp}$.}
    \label{fig:PossibleSituations}
\end{figure}
Similarly, if $(\beta_1-z_L)\cdot\ell^{\perp}<0$, then $X_{T_{\geq2L}-1}^1$ must be within a radius $R$ of $z_L'$. Since $z_L'$ is within distance $R$ from $z_L$, we conclude that $X_{T_{\geq2L}-1}^1$ is within a radius $2R$ of $z_L$, and the event $\mathcal{P}_L$ occurs.

We may therefore assume $\beta_2\cdot\ell^{\perp}$ and $(\beta_1-z_L')\cdot\ell^{\perp}$ are both positive. This situation is depicted on the right in Figure \ref{fig:PossibleSituations}.
If $\alpha_2\cdot\ell^{\perp}>\beta_1\cdot\ell^{\perp}$, then $\cpath^2$ must cross the line $\{x\cdot\ell^{\perp}=\beta_1\cdot\ell^{\perp}\}$ at some point $y$ between $\beta_1$ and $\beta_1+(x_L-2L)\ell$. This crossing point is a distance no more than $\frac{R}{2}$ from $X_n^2$ for some $0\leq n< T_{\leq0}({\bf X}^2)$. Its distance from $\beta_1$ is no more than $R$, and $\beta_1$ is no more than $\frac{R}{2}$ units of distance away from some $X_m^1$ for some $0\leq m\leq T_{\geq 2L}({\bf X}^1)$. Thus, $\mathcal{P}_L$  occurs. 

Finally, assume  $\alpha_2\cdot\ell^{\perp}<\beta_1\cdot\ell^{\perp}$. Then the path taken by $\cpath^1$ from 0 to $\beta_1$ must intersect the path taken by $\cpath^2$ from $\alpha_2$ to $\beta_2$, since they are paths connecting different pairs of opposite corners of the quadrilateral $(0,\beta_2,\beta_1,\alpha_2)$. The point of intersection is no more than $\frac{R}{2}$ units of distance away from $X_m^1$ for some $0\leq m\leq T_{\geq 2L}({\bf X}^1)$ and no more than $\frac{R}2$ away from $X_n^2$ for some $0\leq n< T_{\leq0}({\bf X}^2)$. Thus, $\mathcal{P}_L$ occurs. This finishes the justification of claim \ref{claim:WholeEvent}.

Now \eqref{eqn:282}, together with \eqref{eqn:395}, yields
\begin{equation}\label{eqn:305}
    \lim_{L\to\infty}\Pa^{0,z_L}(\mathcal{O}_L\cup
    \mathcal{P}_L)=\Pa^0(T_{<0}=\infty)\Pa^0(T_{>0}=\infty).
\end{equation}

\hspace{-.25in}{\em Step 3: Handling $\mathcal{O}_L\setminus\mathcal{I}_L$}

In this step, we consider the event $\mathcal{O}_L\setminus\mathcal{I}_L$ and and show that its probability is less than $\frac{1}{2}\Pa^0(G_0^{2L})\Pa^{z_L}(G_{x_L}^0)$.

Because this event does not involve the walks intersecting (and thus ``sharing'' part of the environment), its probability is the same as the probabilities of an analogous event where the two walks are run independently in {\em different} environments. And it is therefore bounded above by the probability of a similar event where two walks are run independently in different environments but are allowed to intersect paths. 
To formalize this idea, let $G_0^{2L,+}$ be the subset of $G_0^{2L}$ on which $X_{T_{\geq2L}}\cdot\ell^{\perp}>z_L$. Let $G_{x_L}^{0,+}$ be the subset of $G_{x_L}^0$ on which $X_{T_{\leq0}}\cdot\ell^{\perp}>0$. Define $G_0^{2L,-}$ and $G_{x_L}^{0,-}$ analogously. And define
\begin{align}\notag
    \Pi_L&:=\{(0=X_0,X_1,\ldots,X_{T_{\geq 2L}}):{\bf X}\in G_0^{2L}\},
    \\
    \notag
    \Pi_{L,+}&:=\{(0=X_0,X_1,\ldots,X_{T_{\geq 2L}}):{\bf X}\in G_0^{2L,+}\},
    \\
    \notag
    \Pi_{L,-}&:=\{(0=X_0,X_1,\ldots,X_{T_{\geq 2L}}):{\bf X}\in G_0^{2L,-}\}.
\end{align}
We will abuse notation by using $\pi$ to denote both a path in one of these sets and the set of vertices in that path. Then
\begin{align}\notag
\Pa^{0,z_L}(\mathcal{O}_L\setminus\mathcal{I}_L)
&=
     \sum_{\pi\in\Pi_{L,+}}\Pa^0(\text{${\bf X}$ takes $\pi$})\Pa^{z_L}(G_{x_L}^{0,-}, T_{\leq0}<T_{\pi})
     \\
     \notag
     &\quad\quad\quad\quad
     +
     \sum_{\pi\in\Pi_{L,-}}\Pa^0(\text{${\bf X}$ takes $\pi$})\Pa^{z_L}(G_{x_L}^{0,+}, T_{\leq0}<T_{\pi})
     \\\notag
     &\leq
     \sum_{\pi\in\Pi_{L,+}}\Pa^0(\text{${\bf X}$ takes $\pi$})\Pa^{z_L}(G_{x_L}^{0,-})
     \\
     \notag
     &\quad\quad\quad\quad
     +
     \sum_{\pi\in\Pi_{L,-}}\Pa^0(\text{${\bf X}$ takes $\pi$})\Pa^{z_L}(G_{x_L}^{0,+})
     \\\notag
     &=\Pa^0(G_0^{2L,+})\Pa^{z_L}(G_{x_L}^{0,-})
    +
    \Pa^0(G_0^{2L,-})\Pa^{z_L}(G_{x_L}^{0,+}) 
    \\\notag
    &\leq \frac12\Pa^0(G_0^{2L})\Pa^{z_L}(G_{x_L}^0).
    \\\notag
    &=\frac12\Pa^0(G_0^{2L})\Pa^{0}(G_0^{-x_L}).
\end{align}
The last inequality comes from the median property of $z_L$, the last equality comes from translation invariance.
Now, using \eqref{eqn:243} and the fact that $\mathcal{I}_L\subset\mathcal{P}_L$, we have 
$$
\limsup_{L\to\infty}\Pa^{0,z_L}(\mathcal{O}_L\setminus\mathcal{P}_L)\leq
\limsup_{L\to\infty}\Pa^{0,z_L}(\mathcal{O}_L\setminus\mathcal{I}_L)\leq\frac12\Pa^0(T_{<0}=\infty)\Pa^0(T_{>0}=\infty)
$$
Hence, due to \eqref{eqn:305},
\begin{equation}\notag
    \frac12\Pa^0(T_{<0}=\infty)\Pa^0(T_{>0}=\infty)\leq\liminf_{L\to\infty}\Pa^{0,z_L}(\mathcal{P}_L). 
\end{equation}

Therefore, to prove \eqref{eqn:222}, it suffices to show
\begin{equation}\label{eqn:464}
    \lim_{L\to\infty}\Pa^{0,z_L}(\mathcal{P}_L\setminus\mathcal{I}_L)=0,
\end{equation} 
and
\begin{equation}\label{eqn:470}
    \lim_{L\to\infty}\Pa^{0,z_L}(\mathcal{I}_L)=0.
\end{equation}

\hspace{-.25in}{\em Step 4: Forcing an unlikely meeting.}

In this step, we show that the above quantities can be compared to the probability of a certain ``meeting event,'' which occurs with vanishing probability. To do this, we exhibit a strategy for ${\bf X}^2$ to hit the path taken by ${\bf X}^1$ at a point appropriately far from $z_L$. 

First, assume that 
\begin{equation}\label{eqn:473}
\Pa^{0,z_L}(\mathcal{P}_L'\setminus\mathcal{I}_L)\geq\frac12\Pa^{0,z_L}(\mathcal{P}_L\setminus\mathcal{I}_L),
\end{equation}
where $\mathcal{P}_L'\subset\mathcal{P}_L$ is the event ${\bf X}^1\in G_0^{2L}$, ${\bf X}^2\in G_{x_L}^0$, and $|X_m^1-X_n^2|\leq2R$ for some $0\leq m\leq T_{\geq 2L}({\bf X}^1)$, $0\leq n < T_{\leq0}({\bf X}^2)$ with $X_m^1\cdot\ell\leq L$. (We will handle the remaining case with an essentially symmetric argument.)

Recall that $M$ is such that for any vertex $y$ in a closed disc of radius $2R$ centered at 0, there is a possible path of length no more than $M$.
Now for a given path $\pi$, define the stopping time
$$T_{\pi,L}':=\inf\left\{n\geq0: \parbox{3.5in}{for some $x\in\pi$ with $x\cdot\ell\leq L$, there is a possible path of length $M$ or less from $X_n$ to $x$}\right\}.$$
Notice that $\mathcal{P}_L'$ implies that $T_{\pi,L}'({\bf X}^2)\leq T_{\leq0}({\bf X}^2)$ for $\pi=(X_n^1)_{n=0}^{T_{\geq2L}}$. This is because if $|X_m^1-X_n^2|\leq2R$, then there is a possible path of length no more than $M$ from $X_n^2$ to $X_m^1$. Therefore,
\begin{align}\notag
\Pa^{0,z_L}(\mathcal{P}_L'\setminus\mathcal{I}_L)
&\leq
\sum_{\pi\in\Pi_L}\Pa^{0,z_L}(\text{${\bf X}^1$ takes $\pi$})\Pa^{0,z_L}(T_{\pi,L}'({\bf X}^2)\leq T_{\leq0}({\bf X}^2)<T_{\pi}({\bf X}^2)\wedge T_{>x_L}({\bf X}^2))
\\\notag
&=     \sum_{\pi\in\Pi_L}\Pa^0(\text{${\bf X}$ takes $\pi$})\Pa^{z_L}(T_{\pi,L}'\leq T_{\leq0}<T_{\pi}\wedge T_{>x_L})
\\\label{eqn:526}
&\leq \sum_{\pi\in\Pi_L}\Pa^0(\text{${\bf X}$ takes $\pi$})\Pa^{z_L}(T_{\pi,L}'<T_{\pi})
.
\end{align}
The first inequality comes from independence of transition probabilities at different sites (condition (\hyperlink{cond:C1}{C1}). Notice that this is where the splitting off of the event $\mathcal{I}_L$ becomes important, because it allows us to leverage independence.

Suppose ${\bf X}$ is such that $T_{\pi,L}'<\infty$ for a path $\pi\in\Pi_L$. Then there is a possible path $\gamma=(y_0=X_{T_{\pi,L}'},y_1,y_2,\ldots,y_k)$ with $k<M$, $y_k\in\pi$, and $y_k\cdot\ell\leq L$. This path does not include any vertices in the path $(X_0,X_1,\ldots,X_{T_{\pi,L}'-1})$, because if such a vertex were included, there would be a shorter path from that vertex to $y_k$, violating the infimum part of the definition of $T_{\pi,L}'$. It may include multiple vertices from $\pi$, but taking $j=\inf\{0\leq i\leq k:y_i\in\pi\}$, we may consider the path $\gamma'=(y_0=X_{T_{\pi,L}'},y_1,y_2,\ldots,y_j)$ that intersects neither $(X_0,X_1,\ldots,X_{T_{\pi,L}'-1})$ nor $\pi$, except at the terminating vertex. 
By the strong Markov property, for every $\omega$ we have
\begin{align}\notag
P_{\omega}^{z_L}\left((X_{T_{\pi,L}'},X_{T_{\pi,L}'+1},\ldots,X_{T_{\pi,L}'+j})=\gamma'\given X_0,\ldots,X_{T_{\pi,L}'}\right)
&=P_{\omega}^{y_0}(\text{{\bf X} takes $\gamma'$}).
\end{align}
Taking conditional expectations on both sides and using the fact that
\begin{equation}\label{eqn:558}
\Pa^a(B|A)=\Ea^a[P_{\omega}^a(B|A)|A]
\end{equation}
for all measurable events $A,B\subset(\Z^2)^{\N}$ with $\Pa^a(A)>0$ (see Lemma \ref{lem:Conditioning}), we have almost surely that $T_{\pi,L}'$ implies
\begin{align}\notag
\Pa^{z_L}\left((X_{T_{\pi,L}'},X_{T_{\pi,L}'+1},\ldots,X_{T_{\pi,L}'+j})=\gamma'\given X_0,\ldots,X_{T_{\pi,L}'}\right)
&=\Ea^a[P_{\omega}^{y_0}(\text{{\bf X} takes $\gamma'$})|X_0,\ldots,X_{T_{\pi,L}'}].
\end{align}
Now by  condition (\hyperlink{cond:C1}{C1}) and the definition of $\gamma'$, this gives us
\begin{align}\label{eqn:ArguingAs}
\Pa^{z_L}\left((X_{T_{\pi,L}'},X_{T_{\pi,L}'+1},\ldots,X_{T_{\pi,L}'+j})=\gamma'\given X_0,X_1,\ldots,X_{T_{\pi,L}'}\right)
&=\Pa^{y_0}(\text{{\bf X} takes $\gamma'$})
\\\label{eqn:545}
&\geq\kappa,
\end{align}
where $\kappa>0$ is the minimum annealed probability of any possible path of length less than $M$. A minimum exists because, up to translation invariance, there are only finitely many possible paths of a given length, and it is positive because by definition, all possible paths have positive annealed probability. Now $y_j$ is of a distance at most $R(k-j)\leq Rk<RM$ from $y_k$, and therefore $y_j\cdot\ell\leq L+RM$. 
If, therefore, ${\bf X}^1$ takes $\pi$ for some $\pi\in\Pi_L$, $T_{\pi,L}'({\bf X}^2)<\infty$, and $(X_{T_{\pi,L}'+1}^2,X_{T_{\pi,L}'+2}^2,\ldots,X_{T_{\pi,L}'+j}^2)=\gamma'$, then $({\bf X}^1,{\bf X}^2)\in\mathcal{M}_L$. Thus, since \eqref{eqn:545} is true whenever $T_{\pi,L}'<\infty$, we have
\begin{equation}\label{eqn:550}
\Pa^{z_L}\left(T_{\pi\cap\{x\cdot\ell<L+RM\}}<\infty\given T_{\pi,L}'<\infty\right)\geq\kappa.    
\end{equation}

We next define another event for $({\bf X}^1,{\bf X}^2)$ that involves the walks intersecting, but is not contained in the strip traversal event. Let the {\em meeting event} $\mathcal{M}_L$ be the event that ${\bf X}^1\in G_0^{2L}$, and ${\bf X}^2$ intersects with the path $(X_n^1)_{n=0}^{T_{\geq2L}}$ at a point $y$ with $y\cdot\ell\leq L+RM$. Note that unlike the  event $\mathcal{I}_L$, our event $\mathcal{M}_L$ does not require that ${\bf X}^2$ complete the event $G_{x_L}^0$, nor does it require that the intersection occur at some $X_n^2$ with $n\leq T_{\leq0}({\bf X}^2)$. However, unlike $\mathcal{I}_L$, $\mathcal{M}_L$ imposes the restriction that the intersection must occur on or near the half of the strip $\{0\leq x\cdot\ell\leq 2L\}$ that is closer to 0.

By \eqref{eqn:550}, we have
\begin{align}\notag
    \Pa^{0,z_L}(\mathcal{M}_L)
    &\geq
    \sum_{\pi\in\Pi_L}
    \Pa^0(\text{{\bf X} takes $\pi$})
    \Pa^{z_L}(T_{\pi,L}'<T_{\pi})\kappa
    \\\notag
    &\overset{\eqref{eqn:526}}{\geq}
    \Pa^{0,z_L}(\mathcal{P}_L'\setminus\mathcal{I}_L)\kappa
    \\\label{eqn:509}
    &\overset{\eqref{eqn:473}}{\geq} \frac12\kappa\Pa^{0,z_L}(\mathcal{P}_L\setminus\mathcal{I}_L),
\end{align}

\hspace{-.25in}{\em Step 4: Finishing the argument.}

We will now show that $\Pa^{0,z_L}(\mathcal{M}_L)$ vanishes. Zerner shows in \cite{Zerner2007} that $\mathcal{I}_L$ has vanishing probability. The argument also works for our event $\mathcal{M}_L$. We summarize it here, applying it to $\mathcal{M}_L$.
Fix $\varepsilon>0$, and suppose the intersection occurs at a point $y$. Either $P_{\omega}^y(A_{\ell})<\varepsilon$ or $P_{\omega}^y(A_{\ell})\geq\varepsilon$. In the former case, a walk from 0 passes through $y$ but still has $T_{\geq L}<T_{<0}$. Zerner shows in \cite{Zerner2007} that the probability of this event has limsup bounded above by $\varepsilon$. In the latter case, a walk started from $z_L$ travels a great distance in direction $-\ell$ (here, a distance at least $L-RM$) and still reaches a point where the probability of $A_{\ell}$ is at least $\varepsilon$. The chance of traveling such a distance in direction $-\ell$ but still having ${\bf X}^2\in A_{\ell}$ approaches 0 as $L\to\infty$. On the other hand, if ${\bf X}^2\in A_{-\ell}$, then $P_{\omega}^{X_n}(A_{\ell})$ must approach 0, being a bounded martingale, and so the probability that it is still above $\varepsilon$ after $\frac{L-RM}{R}$ units of time (long enough to travel distance $L-RM$) approaches 0 as $L\to\infty$. One may then take $\varepsilon$ to 0. Hence we may conclude that
\begin{equation}\notag
\lim_{L\to\infty}\Pa^{0,z_L}(\mathcal{M}_L)=0.
\end{equation}
Since \eqref{eqn:509} is true whenever \eqref{eqn:473} is true, we may conclude that $$
\lim_{L\to\infty}\Pa^{0,z_L}(\mathcal{P}_L\setminus\mathcal{I}_L)\mathbbm{1}_{\text{\eqref{eqn:473} holds}}=0.
$$

Now if \eqref{eqn:473} does not hold, then we can make a nearly symmetric argument. We must have 
\begin{equation}\notag
\Pa^{0,z_L}(\mathcal{P}_L''\setminus\mathcal{I}_L)\geq\frac12\Pa^{0,z_L}(\mathcal{P}_L\setminus\mathcal{I}_L),
\end{equation}
where $\mathcal{P}_L''$ is the event
${\bf X}^1\in G_0^{2L}$, ${\bf X}^2\in G_{x_L}^0$, and $|X_m^1-X_n^2|\leq2R$ for some $0\leq m\leq T_{\geq 2L}({\bf X}^1)$, $0\leq n < T_{\leq0}({\bf X}^2)$ with $X_n^2\cdot\ell\geq L-2R$. Define $\mathcal{M}'$ to be the event that ${\bf X}^2\in G_{x_L}^{0}$ and ${\bf X}^1$ intersects with the path $(X_m^2)_{m=0}^{T_{\leq0}}$ at a point $y$ with $y\cdot\ell\geq L-2R-RM$. As in \eqref{eqn:509}, we can argue that $P^{0,z_L}(\mathcal{M}_L')\geq\frac12\kappa P^{0,z_L}(\mathcal{P}_L\setminus\mathcal{I}_L)$, and we can show as before that $P^{0,z_L}(\mathcal{M}_L')$ vanishes in $L$. Therefore,
$$
\lim_{L\to\infty}\Pa^{0,z_L}(\mathcal{P}_L\setminus\mathcal{I}_L)\mathbbm{1}_{\text{\eqref{eqn:473} does not hold}}=0.
$$

It follows that $\lim_{L\to\infty}\Pa^{0,z_L}(\mathcal{P}_L\setminus\mathcal{I}_L)=0$, which is \eqref{eqn:464}. 
To get \eqref{eqn:470}, note that $\mathcal{I}_L\subset\mathcal{M}_L\cup\mathcal{M}_L'$, so its probability must likewise vanish.
\end{proof}

\section{Random Walks in Dirichlet Environments}\label{sec:Dirichlet}

We now turn our attention to random walks in Dirichlet environments (RWDE). Because the proofs will require graphs other than $\Z^d$, we define RWDE on more general graphs.

Let $\Hg=(V,E,w)$ be a weighted directed graph with vertex set $V$, edge set
$~E\subseteq V\times V$, and a positive-valued weight function $w:E\to\R$. 
To the weighted directed graph $\Hg$, we can associate the Dirichlet measure $\Pm_{\Hg}$ on $(\Omega_{V},\mathcal{F}_{V})$, which we now describe. 

Recall the definition of the Dirichlet distribution: for a finite set $I$, take parameters $\alpha =(\alpha_i)_{i\in I}$, with $\alpha_i>0$ for all $i$. The Dirichlet distribution with these parameters is a probability distribution on the simplex $\Delta_I:=\{(x_i)_{i\in I}:\sum_{i\in I}x_i=1\}$ with density
\begin{equation*}D\left((x_i)_{i\in I}\right)=C(\alpha )\prod_{i\in I}x_i^{\alpha_i-1},\end{equation*}
where $C(\alpha )$ is a normalizing constant.

Define $\Pm_{\Hg}$ to be the measure on $\Omega_V$ under which transition probabilities at the various vertices $x\in V$ are independent, and for each vertex $x\in V$, $(\omega(x,y))_{(x,y)\in E}$ is distributed according to a Dirichlet distribution with parameters $(w(x,y))_{(x,y)\in E}$. With $\Pm_{\Hg}$-probability 1, $\omega(x,y)>0$ if and only if $(x,y)\in E$ for all $x,y\in V$. We will call a random environment chosen according to $\Pm_{\Hg}$ a {\em Dirichlet environment on $\Hg$}.  We will use $E_{\Hg}$ to denote the associated expectation, and $\Pa_{\Hg}^x$ and $\Ea_{\Hg}^x$ to denote the annealed measure and expectations.

Now let $\mathscr{N}\subset\Z^d$ be a finite set such that $\sum_{i=1}^{\infty}(\mathscr{N}\cup\{0\})=\Z^d$.
Let $\Gr$ be a weighted directed graph with vertex set $\Z^d$, and let $(\alpha_y)_{y\in\mathscr{N}}$ be positive weights. Let $\Gr=(\Z^d,E,w)$ be the weighted directed graph with vertex set $\Z^d$, edge set $E:=\{(x,y)\in\Z^d\times\Z^d:y-x\in\mathscr{N}\}$, and weight function $w$ with $w(x,y)=\alpha_{y-x}$ for all $(x,y)\in E$ (the condition on $\mathscr{N}$ simply ensures that $\Gr$ is strongly connected). Then $\Pm_{\Gr}$ is the law of a Dirichlet environment on $\Z^d$ satisfying (\hyperlink{cond:C1}{C1}), (\hyperlink{cond:C2}{C2}), and (\hyperlink{cond:C3}{C3}), and $\Pa_{\Gr}^0$ is the corresponding annealed measure for a walk started at 0. 

It is known in nearest-neighbor RWDE that for a given direction $\ell\in S^{d-1}$, transience and recurrence in direction $\ell$ under $\Pa_{\Gr}^0$ are characterized by the relationship between $\ell$ and the annealed drift. 

\begin{thm*}[{\cite[Theorem 1]{Sabot&Tournier2016}}]
Let $\Pa_{\Gr}^0$ be the measure of a nearest-neighbor RWDE on $\Z^d$. Let $\Delta=\Ea^0[X_1]$ be the annealed drift, and let $\ell\in S^{d-1}$. Then $\Pa_{\Gr}^0(A_{\ell})=1$ if and only if $\ell\cdot\Delta>0$; otherwise, $\Pa_{\Gr}^0(A_{\ell})=0$. 
\end{thm*}

Our goal is to extend this theorem to RWDE with bounded jumps.

\begin{thm}\label{thm:DirichletFullStatement}
Let $\Pa_{\Gr}^0$ be the measure of a RWDE with bounded jumps on $\Z^d$. Let $\Delta=\Ea^0[X_1]$ be the annealed drift, and let $\ell\in S^{d-1}$. Then $\Pa_{\Gr}^0(A_{\ell})=1$ if and only if $\ell\cdot\Delta>0$; otherwise, $\Pa_{\Gr}^0(A_{\ell})=0$. 
\end{thm}

As Tournier points out in \cite{Tournier2015}, many of the arguments used in the proof of the theorem from \cite{Sabot&Tournier2016} do not rely on the nearest-neighbor assumption, and therefore already work for RWDE with bounded jumps as well. In particular, Theorem \ref{thm:DirichletFullStatement} is known to be true provided $\Delta\neq0$ and $d\neq2$. 

If $\Delta\neq0$ and $d=2$, we know from \cite{Tournier2015} that $\ell\cdot\Delta>0$ implies $\Pa_{\Gr}^0(A_{\ell})>0$, and from \cite[Theorem 1.8]{Drewitz&Ramirez2010} that $\ell\cdot\Delta=0$ implies $\Pa_{\Gr}^0(A_{\ell})=0$ (the arguments in \cite{Drewitz&Ramirez2010} are given for the nearest-neighbor case, but can be easily modified to work for our bounded-jump model). From here, our 0-1 law of Theorem \ref{thm:01Law} allows us to reach the conclusion of Theorem \ref{thm:DirichletFullStatement}. 

The only remaining case is where $\Delta=0$. In the nearest-neighbor case, $\Delta=0$ implies a symmetry that forces $\Pa_{\Gr}^0(A_{\ell})=\Pa_{\Gr}^0(A_{-\ell})$ for all directions $\ell$. The 0-1 laws of \cite{Zerner&Merkl2001} for $d=2$ and of \cite{Bouchet2013} for $d\geq3$ then yield the conclusion $\Pa_{\Gr}^0(A_{\ell})=0$ for all $\ell$. 
In the bounded-jumps case, zero drift does not imply symmetry, so even the 0-1 law of Theorem \ref{thm:01Law} is not by itself enough to prove the theorem. Theorem \ref{thm:DirichletFullStatement} will be proven if we can prove the following theorem, which will rely on Theorem \ref{thm:01Law} for the case $d=2$. 

\begin{thm}\label{thm:Dirichlet}
Let $\Pa_{\Gr}^0$ be the measure of a RWDE with bounded jumps on $\Z^d$. Let $\Delta=0$, and let $\ell\in S^{d-1}$. Then $\Pa_{\Gr}^0(A_{\ell})=0$. 
\end{thm}

As with many proofs of results in RWDE, our proof involves comparing the graph $\Gr$ to a sequence of larger and larger finite graphs ($\Hg_{N,L}$), which look like $\Gr$ except possibly near boundaries, and applying a key lemma from \cite{Sabot2011} (we only state the part we need).
Let $\Hg=(V,E,w)$ be a weighted directed graph. For a site $x\in V$, define the \textit{divergence} of $x$
by
$\text{div}(x)=\sum_{(x,y)\in E}w(x,y)-\sum_{(y,x)\in E}w(y,x)$. Say $x$ has \textit{zero divergence} if $\text{div}(x)=0$, and say $\Hg$ is \textit{divergence-free} if every vertex in $V$ has zero divergence.
\begin{lem}[See \cite{Sabot2011}, Lemma 1, \cite{Sabot&Tournier2011}, Lemma 1]\label{cor:loopreversal}
Let $\Hg=(V,E,w)$ be a divergence-free weighted directed graph, and let $x,y\in V$ such that there is an edge $e$ from $y$ to $x$ in $\Hg$. Then, letting $\tilde{T}_x$ denote the first positive hitting time of $x$, $\Pa_{\Hg}^x(X_{\tilde{T}_x-1}=y)=\frac{w(y,x)}{\sum_{v\in V}w(v,x)}$.
\end{lem}
The finite graphs we construct are closely related to those constructed by Tournier in \cite{Tournier2015} for the characterization of transience in the nonzero-drift case, and we invoke a key property of these graphs that is proven in \cite{Tournier2015} and also applies to our graphs. However, the graphs from \cite{Tournier2015} are slightly altered to suit our argument, which is new, though of a very similar flavor. 

A very significant difference between the argument in \cite{Tournier2015} and ours is that the former need only be given for directions with rational slopes, and the result follows for arbitrary directions as an immediate consequence.
That is, let $S_r^{d-1}:=\left\{\frac{\vecu}{|\vecu|}:\vecu\in\Z^d\setminus\{0\}\right\}\subset S^{d-1}$ be the set of vectors in the unit sphere $S^{d-1}$ that have all rational slopes. The finite graphs used in 
\cite{Tournier2015} are naturally set up to work for directions in $S_r^{d-1}$, but once the result is proven, extending to arbitrary directions follows immediately from convexity of the set of transient directions. For us, extending the result from rational directions to all $\ell\in S^{d-1}$ is not so immediate. Indeed, the following conjecture remains open for general RWRE.
\begin{conj}\label{conj1}
Let $\Pa^0$ be the law of an i.i.d. RWRE on $\Z^d$, and let $S^{d-1}$ be the set of a unit vectors in $\R^d$. For $\ell\in S^{d-1}$, let $A_{\ell}^0$ be the event that $\lim_{n\to\infty}X_n\cdot\ell=\infty$, but there is no neighborhood $U\in S^{d-1}$ containing $\ell$ such that for all $\ell'\in U$, $\lim_{n\to\infty}X_n\cdot\ell'=\infty$. Then for all $\ell\in S^{d-1}$, $\Pa^0(A_{\ell}^0)=0$.
\end{conj}
In the nearest-neighbor case of RWDE, Conjecture \ref{conj1} is seen to be true from \cite[Theorem 1]{Sabot&Tournier2016}, and in the bounded-jump case it will follow from Theorem \ref{thm:DirichletFullStatement}, once it is proven (again, it only remains to prove Theorem \ref{thm:Dirichlet}). 
However, because we cannot rely on the truth of Conjecture \ref{conj1}, proving Theorem \ref{thm:Dirichlet} for all directions $\ell\in S_r^{d-1}$ is not sufficient to prove it for all directions $\ell\in S^{d-1}$, even though $S_r^{d-1}$ is dense in $S^{d-1}$. For $\ell\in S^{d-1}\setminus S_r^{d-1}$, we must rule out the possibility that  a walk could with positive probability be transient in direction $\ell$ while recurrent in all directions not parallel to $\ell$.

For the sake of readability, we will first prove Theorem \ref{thm:Dirichlet} for rational slopes. However, we do not know of a way to generalize directly to arbitrary directions. 
Rather, the generalization will require going through the same argument more carefully, choosing directions $\vecv\in S_r^{d-1}$ sufficiently close to $\ell$ to satisfy certain properties and constructing graphs in terms of these $\vecv$. We will describe necessary differences as they come up. This approach should be easier to follow, as the ideas involved in constructing the graphs and leveraging the time reversal lemma are quite separate from the ideas involved in comparing arbitrary directions with directions in $S_r^{d-1}$.
 
\subsection{Preliminaries}
 Because we are discussing multiple directions, we must replace our notation $T_{\leq a}$ for hitting times of half-spaces with the slightly more cumbersome
$$
T_{\leq a}^{\ell}=T_{\leq a}^{\ell}({\bf X}):=\inf\{n\geq0:(X_n\cdot\ell)\leq a\},
$$
and similarly for $<$, $\geq$, and $>$. We use this notation even for the proof that assumes $\ell\in S_r^{d-1}$ in order to facilitate comparisons later. 

Moreover, for $\ell\in S^{d-1}$ we will also need to define the ``lateral hitting times''
$$H_{\geq a}^{\ell}:=\inf\{n\geq0:X_n\cdot\ell^{\perp}\geq a\text{ for some }\ell\in S^{d-1}\text{ with }\ell^{\perp}\perp\ell\}.$$

For a vertex set $V$ and $v\in V$, define
$$
T_v=T_v({\bf X}):=\inf\{n\geq0:X_n=v\}.
$$

Finally, for any stopping time defined as the first $n\geq0$ satisfying a certain condition, we use the same notation but with a tilde ($\sim$) over it to denote the corresponding positive stopping time: that is, the first $n>0$ satisfying the same condition. 

We will use the following lemma, which is also an ingredient in Kalikow's 0-1 law. It is proven in Appendix \ref{app:Kalikow}.

\begin{lem}\label{lem:NoGettingStuckInaStrip}
    Let $\Pa^0$ be the annealed measure of a RWRE on $\Z^d$ satisfying assumptions (\hyperlink{cond:C1}{C1}), (\hyperlink{cond:C2}{C2}), and (\hyperlink{cond:C3}{C3}). Then for every $\ell\in S^{d-1}$ and $a<b\in\R$,
    \begin{equation}\label{eqn:11}
    \Pa^0(\#\{n\geq0: X_n\cdot\ell\geq a\}=\infty, T_{\geq b}^{\ell}=\infty)=0.
    \end{equation}
\end{lem}

\subsection{Rational slopes}
We now state and prove Theorem \ref{thm:Dirichlet} for directions with rational slopes.

\begin{thm}\label{thm:DirichletRational}
Let $\Pa_{\Gr}^0$ be the measure of a RWDE with bounded jumps on $\Z^d$. Let $\Delta=0$, and let $\ell\in S_r^{d-1}$. Then $\Pa_{\Gr}^0(A_{\ell})=0$. 
\end{thm}

\begin{proof}[Proof of Theorem \ref{thm:DirichletRational}] 
Let $\ell\in S_r^{d-1}$. Assume for a contradiction that $\Pa_{\Gr}^0(A_{-\ell})>0$. Then, as in \cite[page 765]{Kalikow1981}, we have $\Pa_{\Gr}^0(T_{>0}^{\ell}=\infty)>0$, from which it easily follows that $\alpha:=\Pa_{\Gr}^0(\tilde{T}_{\geq0}^{\ell}=\infty)>0$. 
By Lemma \ref{lem:NoGettingStuckInaStrip},
$\Pa_{\Gr}^0(T_{\leq-L}^{\ell}=\tilde{T}_{\geq0}^{\ell}=\infty)=0$, so $\Pa_{\Gr}^0(T_{\leq-L}^{\ell}<\tilde{T}_{\geq0}^{\ell}=\infty)=\alpha$. For $P_{\Gr}$-almost every environment, it is possible, with positive probability, for a walk to hit the half-space $\{x\cdot\ell\leq-L\}$ and then return to $\{x\cdot\ell\geq0\}$. Therefore,
we have $\Pa_{\Gr}^0(T_{\leq-L}^{\ell}<\tilde{T}_{\geq0}^{\ell})>\alpha$ for all $L\geq0$. And on the event $\{T_{\leq-L}^{\ell}<\tilde{T}_{\geq0}^{\ell}\}$, there is necessarily some $a$ such that $T_{\leq-L}^{\ell}<\tilde{T}_{\geq0}^{\ell}\wedge H_{\geq a}^{\ell}$. It therefore follows that for any $L>0$, there exists $K=K(L)>0$ such that 
\begin{equation}\label{eqn:585}
\Pa_{\Gr}^0(T_{\leq-L}^{\ell}<\tilde{T}_{\geq0}^{\ell}\wedge H_{\geq \frac12K}^{\ell})>\alpha.
\end{equation}

For $L\geq R$, let $K=K(L)$ be an increasing function satisfying \eqref{eqn:585} for all L. Let $\vecu$ be a constant multiple of $\ell$ such that $\vecu\in\Z^d$. Then let $(\vecuk,\vecuk_2,\ldots,\vecuk_d)$ be an orthogonal basis for $\R^d$ such that $\vecuk_i\in\Z^d$ for all $i$. Let $N$ be large enough that $N|\vecu_i|\geq K$ for all $i$.

We will define a graph $\Hg_{N,L}$ in nearly the same way as the $G_{N,L}$ defined by Tournier in \cite{Tournier2015}. Consider the cylinder
$$
C_{N,L}:=\{x\in\Z^d:0\leq x\cdot\ell\leq L\}/(N\Z\vecu_2+\cdots+N\Z\vecu_d).
$$
This is the slab $\mathcal{S}_{N,L}:=\{0\leq x\cdot\ell\leq L\}\cap\Z^d$ where vertices that differ by $N\vecu_i$ for some $i\in\{2,\ldots,d\}$ are identified. We note that \cite{Tournier2015} uses $L|\vecu|$ rather than $L$ here. We use $L$ for reasons related to our plans for generalizing the proof to $\ell\notin S_r^{d-1}$. 

Now define the graph $\Hg_{N,L}$ with vertex set 
$$V_{N,L}:=C_{N,L}\cup\{M\}\cup\{\partial\},$$ 
where $M$ and $\partial$ are new vertices (in \cite{Tournier2015}, $R$ and $\partial$ are used, but in this paper, $R$ refers to the jump range, as defined in (\hyperlink{cond:C3}{C3}), and edges of $\Hg_{N,L}$ are of the following types:
\begin{enumerate}[1.]
    \item edges induced by those of $\Gr$ inside $C_{N,L}$;
    \item If $x\in C_{N,L}$ corresponds to a vertex $x'\in\mathcal{S}_{N,L}$, there is 
    \begin{enumerate}[(a)]
        \item an edge from $x$ to $\partial$ for each $y\in\mathscr{N}$ such that $(x'+y)\cdot\ell<0$,
        \item an edge from $\partial$ to $x$ for each $y\in-\mathscr{N}$ such that $(x'+y)\cdot\ell<0$,
        \item an edge from $x$ to $M$ for each $y\in\mathscr{N}$ such that $(x'+y)\cdot\ell>L$,
        \item an edge from $M$ to $x$ for each $y\in-\mathscr{N}$ such that $(x'+y)\cdot\ell>L$,
    \end{enumerate}
    \item A new ``special'' edge from $M$ to $\partial$ and one from $\partial$ to $M$.
\end{enumerate}
Weights of all edges but the last two are induced by the corresponding weights in $\Gr$. Note that several edges may share the same head and tail. If that is the case, identify such edges into one edge whose weight is the sum of all of the original weights in order to create a graph that is not a multigraph and fits our definitions (there is also a way to define RWDE on a multigraph by keeping track of vertices visited and edges taken, and if we used such a definition, the identification of multiple edges would not affect the distribution of the vertex path). By construction,  all vertices in $C_{N,L}$ have zero divergence. It remains to describe the weights of the new edges connecting $\partial$ and $M$ (the paper \cite{Tournier2015} only defines an edge from $M$ to $\partial$). It is shown in \cite[p. 722]{Tournier2015} that the quantity
$$
\left(\sum\text{weights of edges in 2(c)}\right)
-
\left(\sum\text{weights of edges in 2(a)}\right)
$$
is a positive multiple of the dot product of $v$ with the annealed drift. Thus, because of our assumption the annealed drift is zero, the two sums are equal. Note that by the shift-invariant structure of the graph $\Gr$, the sum on the left is also the weight exiting $\partial$ by edges in 2(b). Similarly, the sum on the right is also the weight exiting $M$ by edges in 2(d). Hence the total weights of edges in 2(a), 2(b), 2(c), and 2(d) are all the same. Because weights in 2(a) and 2(b) are the same, $\partial$ has zero divergence, and because 2(c) and 2(d) are the same, $M$  has zero divergence. In order to preserve the divergence-free character of the graph, we give both of the special edges the same weight $W$, which we take to be the value of each of the two sums above. It follows from well known properties of Dirichlet random variables that when the walk is started at either of the endpoints, its first step is along the special edge to the other endpoint with annealed probability $\frac12$. Figure \ref{fig:graphHM} shows an example of the graph $\Hg_{N,L}$. (Because figure \ref{fig:graphHM} is also intended to be used for the argument for Theorem \ref{thm:Dirichlet}, it uses $\vecv$ rather than $\ell$ in its labeling. For the purpose of the current argument, simply take $v=\ell$.)

\begin{figure}
    \centering
\includegraphics[height=3.5in]{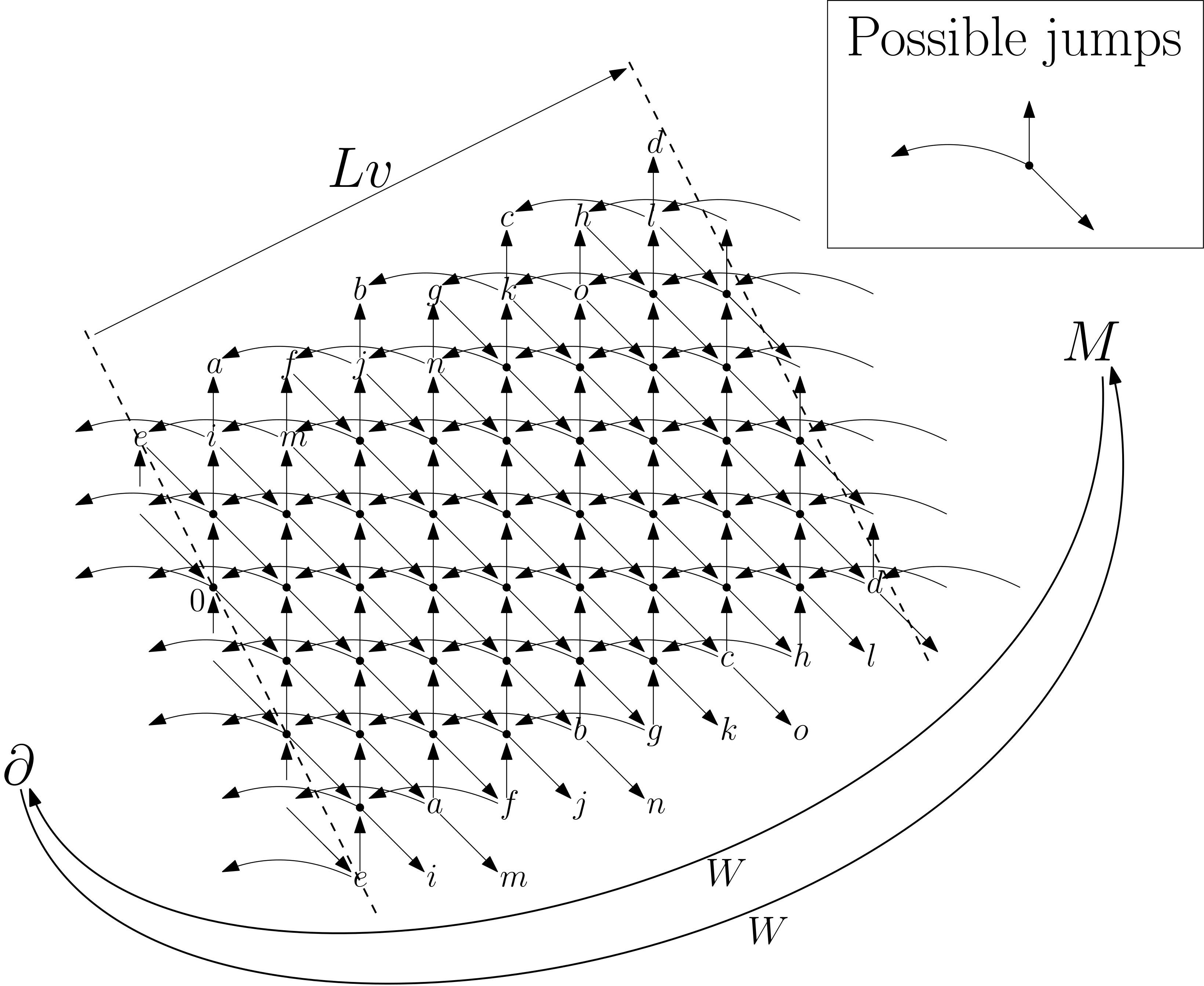}
    \caption{Graph $\Hg_{N,L}$. Here $\mathscr{N}=\{(0,1),(1,-1),(-2,0)\}$, and $v=(2,1)$. Boundary conditions in direction perpendicular to $v$ are periodic; vertices labeled with the same letters are identified. Arrows to and from the main part of the graph on the left are understood to originate from or terminate at $\partial$, and similarly with $M$ on the right side.}
    \label{fig:graphHM}
\end{figure}

Define the stopping time $\tau=\inf\{n\in\N:X_n=\partial, X_{n-1}=M\}$. Note that $\{\tilde{T}_{\partial}=\tau\}=\{X_{\tilde{T}_{\partial}-1}=M\}$ is the event that the first return to zero is by the special edge. 

We note, by Lemma \ref{cor:loopreversal},
\begin{equation*}\Pa_{\Hg_{N,L}}^{\partial}(\tilde{T}_{\partial}= \tau)=\frac{1}{2}
\end{equation*}
On the other hand, we also have $\Pa_{\Hg_{N,L}}^{\partial}(X_1=M)=\frac{1}{2}$.
Now by considering the possibility that the first step from $\partial$ is to $M$ (by the special edge) and the possibility that the first step from $\partial$ is not to $M$, we get
\begin{align}\notag\Pa_{\Hg_{N,L}}^{\partial}(\tilde{T}_{\partial}= \tau)&\leq
\Pa_{\Hg_{N,L}}^{\partial}(X_1=M,~\tilde{T}_{\partial}= \tau)+\Pa_{\Hg_{N,L}}^{\partial}(X_1\neq  M,T_{M}<\tilde{T}_{\partial})
\\\label{eqn:783}
&=\Pa_{\Hg_{N,L}}^{\partial}(X_1=M)\Pa_{\Hg_{N,L}}^{M}(T_{\partial}= \tau)+\Pa_{\Hg_{N,L}}^{\partial}(X_1\neq  M,T_{M}<\tilde{T}_{\partial}).
\end{align}
The equality comes from the Markovian property of the  quenched measure, Lemma \ref{lem:Conditioning}, and independence of sites, arguing as we did for \eqref{eqn:ArguingAs}.
Now \eqref{eqn:783} can be rewritten as
\begin{equation}\label{eqn:804}
\frac{1}{2}
\leq 
\frac{1}{2}\Pa_{\Hg_{N,L}}^{M}(\tilde{T}_{\partial}= \tau)
+
\Pa_{\Hg_{N,L}}^{\partial}(X_1\neq  M,T_{M}<\tilde{T}_{\partial}),
\end{equation}

We claim that the term $\Pa_{\Hg_{N,L}}^{\partial}(X_1\neq  M,T_{M}<\tilde{T}_{\partial})$ approaches 0 as $L$ and $K$ increase. Let $B=B(L,K)$ be a box of radius $\frac{L\wedge K}{3}$ around 0, and for $x\in C_{N,L}$, let $x+B$ be the set of vertices in $C_{N,L}$ that can be written as $x+y$ for some $y\in B$. Note that for $x\in C_{N,L}$, the dot product with $\ell$ is well defined, since vertices in $\mathcal{S}_{N,L}$ that are identified to form $C_{N,L}$ have the same dot product with $\ell$.
Then for sufficiently large $L$,
\begin{align}
\notag
    \Pa_{\Hg_{N,L}}^{\partial}(X_1\neq  M,T_{M}<\tilde{T}_{\partial})
    &=\sum_{\substack{x\in C_{N,L},\\0\leq x\cdot\ell\leq R}}\Pa_{\Hg_{N,L}}^{\partial}(X_1=x)\Pa_{\Hg_{N,L}}^x(T_M<T_{\partial})
    \\\notag
    &\leq\sum_{\substack{x\in C_{N,L},\\0\leq x\cdot\ell\leq R}}\Pa_{\Hg_{N,L}}^{\partial}(X_1=x)\Pa_{\Hg_{N,L}}^x(T_{(x+B)^c}<T_{\partial})
    \\\notag
    &\leq\sum_{\substack{x\in C_{N,L},\\0\leq x\cdot\ell\leq R}}\Pa_{\Hg_{N,L}}^{\partial}(X_1=x)\Pa_{\Gr}^0(T_{B^c}<T_{\leq-R}^{\ell})
    \\\label{eqn:669}
    &=\Pa_{\Gr}^0(T_{B^c}<T_{\leq-R}^{\ell})
\end{align}
The first equality comes from the strong Markov property and independence of sites. The first inequality holds as long as $L$ is large enough that $M\notin x+B$. To get the second inequality, note that a finite path from $x$ that stays in $x+B$ until the last step does not use the periodic boundary conditions (provided $\frac{L\wedge K}{3}>R$), and so it has the same probability as a corresponding path in $\Gr$. And for $x\in C_{N,L}$ with $x\cdot\ell\leq R$,
a walk from $x$ on $\Hg_{N,L}$ that leaves $x+B$ without hitting $\partial$ corresponds to a walk on $\Gr$ (which we may take to start at 0 by translation invariance) that leaves $B$ without traveling $x\cdot\ell$ or more units (of distance in $\R^d$) in direction $-\ell$. Since $x\cdot\ell\leq R$ for all $x$ with $\Pa_{\mathcal{H}_{N,L}}^{\partial}(X_1=x)>0$, the second inequality follows. The final equality comes from pulling the second term out of the sum, which is then equal to 1.

To prove our claim, we must show that \eqref{eqn:669} goes to 0 as $L$ increases (along with $K$). Let $\varepsilon>0$. By assumption, $\Pa_{\Gr}^0(A_{-\ell})>0$. By Theorem \ref{thm:01Law} for $d=2$, or by the 0-1 law of Bouchet for $d\geq3$ in \cite{Bouchet2013} (where, as Tourner points out in \cite{Tournier2015}, the proof works for bounded jumps), this means $\Pa_{\Gr}^0(A_{-\ell})=1$. Thus, $\Pa_{\Gr}^0(T_{\leq -R}^{\ell}<\infty)=1$. Now take an increasing sequence $(Q_r)$ of finite sets converging to $\Z^d$. Then the event $\{T_{\leq-R}^{\ell}<\infty\}$ is the limit as $r$ increases (i.e., the union over all $r$) of the events $\{T_{\leq-R}^{\ell}<T_{Q_r^c}\}$. Let $Q=Q(\varepsilon)$ be one such $Q_r$ large enough that $\Pa_{\Gr}^0(T_{Q^c}<T_{\leq-R}^{\ell})<\varepsilon$. Note that although $Q$ depends on $\varepsilon$, it does not depend on $L$. Thus, for large enough $L$, $B$ contains $Q$, so that
\begin{equation}\label{eqn:874}
\{T_{B^c}\leq T_{\leq-R}^{\ell}\}\subset\{T_{Q^c}\leq T_{\leq-R}^{\ell}\}.
\end{equation}
It follows that, for large enough $L$,
$$
\Pa_{\Gr}^0(T_{B^c}\leq T_{\leq-R}^{\ell})\leq \Pa_{\Gr}^0(T_{Q^c}\leq T_{\leq-R}^{\ell})<\varepsilon.
$$
Since this can be true for arbitrary $\varepsilon>0$, the right side of \eqref{eqn:669} goes to 0, and therefore so does $\Pa_{\Hg_{N,L}}^{\partial}(X_1\neq  M,T_{M}<\tilde{T}_{\partial})$.

Next, we will show that $\Pa_{\Hg_{N,L}}^M(T_{\partial}=\tau)$ is bounded away from 1 as $M$ increases. We have

\begin{align*}
    \Pa_{\Hg_{N,L}}^M(T_{\partial}\neq \tau)
    &\geq \sum_{x\in C_{N,L},L-R\leq x\cdot\ell\leq L}\Pa_{\Hg_{N,L}}^M(X_1=x)\Pa_{\Hg_{N,L}}^{x}(T_{\partial}<T_{>x\cdot\ell}^{\ell})
    \\
    &\geq \sum_{x\in C_{N,L},L-R\leq x\cdot\ell\leq L}\Pa_{\Hg_{N,L}}^M(X_1=x)\Pa_{\Gr}^0(T_{\leq-L}^{\ell}<\tilde{T}_{\geq0}^{\ell}\wedge H_{\geq \frac12K}^{\ell})
    \\
    &>\sum_{x\in C_{N,L},L-R\leq x\cdot\ell\leq L}\Pa_{\Hg_{N,L}}^M(X_1=x)\alpha
    \\
    &=\frac12\alpha.
\end{align*}
The first inequality comes from the strong Markov property and independence of sites. To get the second inequality, note that the probability $\Pa_{\Hg_{N,L}}^{x}(T_{\partial}<T_{>x\cdot\ell})$ is greater than the probability, starting from $x$, that a walk on $\Hg_{N,L}$ reaches $\partial$ without ever traveling more than $\frac{N}{3}$ units in any direction perpendicular to $\vecuk$. Since this event precludes the walk from using the periodic boundary conditions, (and because weights to $\partial$ in $\Hg_{N,L}$ are the same as the weights from corresponding sites to the set $\{y:y\cdot\vecuk<0\}$) its probability is the same as the probability that a walk in $\Gr$ travels more than $x\cdot\ell$ units in direction $-\vecuk$ without ever traveling more than $\frac{N}{3}$ units in any perpendicular direction. Since $x\cdot\vecuk\leq L$, the second inequality follows. The third inequality comes from \eqref{eqn:585}, and the equality comes from the expectation of a beta random variable. 

Now taking the limsup in \eqref{eqn:804} as $M\to\infty$ yields the contradiction

$$
\frac{1}{2}
\leq 
\frac{1}{2}\left(1-\frac12\alpha\right)
<
\frac{1}{2}.
$$
\end{proof}

\subsection{Generalizing to directions in $S^{d-1}\setminus S_r^{d-1}$}

We now describe how to prove Theorem \ref{thm:DirichletRational} for directions that do not necessarily have rational slopes. 

The graph constructed in \cite{Tournier2015} is used to analyze a direction $\ell$ with rational slopes, and uses the rationality in a significant way. Rather than attempt to construct and analyze an analogous graph for an irrational direction $\ell\in S^{d-1}$, we use a sequence of rational slopes $\vecvk\in S_r^{d-1}$ approaching $\ell$. The following lemma is simple, but important.

\begin{lem}\label{lem:598}
Fix $\ell\in S^{d-1}$, $h>0$, and $L'>L>0$. For $\vecv$ close enough to $\ell$, any $x\in\R^d$ with $x\cdot\ell\geq L'$ and $x\cdot\vecv\leq L$ must necessarily have $x\cdot\ell^{\perp}>h$ for some $\ell^{\perp}\perp\ell$. 
\end{lem}
\begin{proof}
Choose a unit vector $\vecv$ close to $\ell$ and let $\ell^{\perp}\in S^{d-1}$ be the unit vector perpendicular to $\ell$ such that $\vecv=a\ell-\sqrt{1-a^2}\ell^{\perp}$, where $a=\vecv\cdot\ell$. Then $a\upto1$ as $\vecv\to\ell$, and $\ell-\vecv=(1-a)\ell+\sqrt{1-a^2}\ell^{\perp}$. By writing $x\cdot(\ell-\vecv)$ in different ways, we get
$$
(1-a) x\cdot\ell + \sqrt{1-a^2} x\cdot\ell^{\perp} = x\cdot\ell-x\cdot\vecv.
$$
From this we get
\begin{align*}
    \sqrt{1-a^2}x\cdot\ell^{\perp}&=ax\cdot\ell-x\cdot\vecv
    \\
    &\geq aL'-L.
\end{align*}
For $\vecv$ sufficiently close to $\ell$, $a$ is close enough to 1 that this gives us
$$
\sqrt{1-a^2}x\cdot\ell^{\perp}\geq\frac12(L'-L)
$$
and
$$
x\cdot\ell^{\perp}\geq\frac{1}{2\sqrt{1-a^2}}(L'-L).
$$
Taking $\vecv$ close to $\ell$ makes $a$ close to 1, which suffices to prove the lemma. 
\end{proof}

We now proceed with the proof, describing only the parts where it differs from the proof of Theorem \ref{thm:DirichletRational}. 
\begin{proof}[Proof of Theorem \ref{thm:Dirichlet}]

The first challenge is to get the same bound as in \eqref{eqn:585}, but for a direction $\vecv$ with rational slopes. We will show that for any $L$, there is a unit vector $\vecv=\vecv(L)\in S_r^{d-1}$ close enough to $\ell$ and a $K=K(L)$ large enough that 
\begin{equation}\label{eqn:594}
\Pa_{\Gr}^0(T_{\leq-L}^{\vecv}<\tilde{T}_{\geq0}^{\vecv}\wedge H_{\geq \frac12K}^{\vecv})>\alpha.
\end{equation}
Fix $L>0$, and choose any $L'>L$. Let $K'$ be such that 
\begin{equation}\notag
\Pa_{\Gr}^0(T_{\leq-L'}^{\ell}<\tilde{T}_{\geq0}^{\ell}\wedge H_{\geq \frac12K'}^{\ell})>\alpha.    
\end{equation}
(Such a $K'$ exists by \eqref{eqn:585}.) 
Now on the event $\{T_{\leq-L'}^{\ell}<\tilde{T}_{\geq0}^{\ell}\}$, there is necessarily an open neighborhood around $\ell$ such that for any $\vecv$ in the neighborhood, $T_{\leq-L'}^{\ell}<\tilde{T}_{\geq0}^{\vecv}$. This is because the walk only hits finitely many points before $T_{\leq-L'}$, and each such point $x$ (other than 0) has $x\cdot\ell<0$, so that for $\vecv$ close enough to $\ell$, $x\cdot\vecv<0$. Hence
$$
\lim_{\vecv\to\ell}\Pa_{\Gr}^0(T_{\leq-L'}^{\ell}<\tilde{T}_{\geq0}^{\ell}\wedge \tilde{T}_{\geq0}^{\vecv}\wedge H_{\geq \frac12K'}^{\ell})=\Pa_{\Gr}^0(T_{\leq-L'}^{\ell}<\tilde{T}_{\geq0}^{\ell}\wedge H_{\geq \frac12K'}^{\ell}).
$$
In particular, for $\vecv$ close enough to $\ell$, 
\begin{equation}\label{eqn:591}
\Pa_{\Gr}^0(T_{\leq-L'}^{\ell}<\tilde{T}_{\geq0}^{\ell}\wedge \tilde{T}_{\geq0}^{\vecv}\wedge H_{\geq \frac12K'}^{\ell})>\alpha.
\end{equation}
Now let $\vecv\in S^{d-1}$ have rational slopes, satisfy \eqref{eqn:591}, and also be close enough to $\ell$ that if $x\cdot\ell\geq L'$ and $x\cdot\vecvk\leq L$, then $x\cdot\ell^{\perp}\geq K'$ for some $\ell^{\perp}\perp\ell$ (this is possible by Lemma \ref{lem:598}).
Choose $K$ large enough that any $y$ with $y\cdot\vecvk^{\perp}\geq \frac{K}{2}$ for any $\vecvk^{\perp}\perp\vecvk$ is necessarily outside the set 
$$Z:=\left\{-L'-R\leq x\cdot\ell\leq0, x\cdot\vecvk\leq0, x\cdot\ell^{\perp}\leq\frac{K'}{2}\text{ for all }\ell^{\perp}\perp\ell\right\}.$$
\renewcommand{\floatpagefraction}{.8}
\begin{figure}
    \centering
\includegraphics[height=4in]{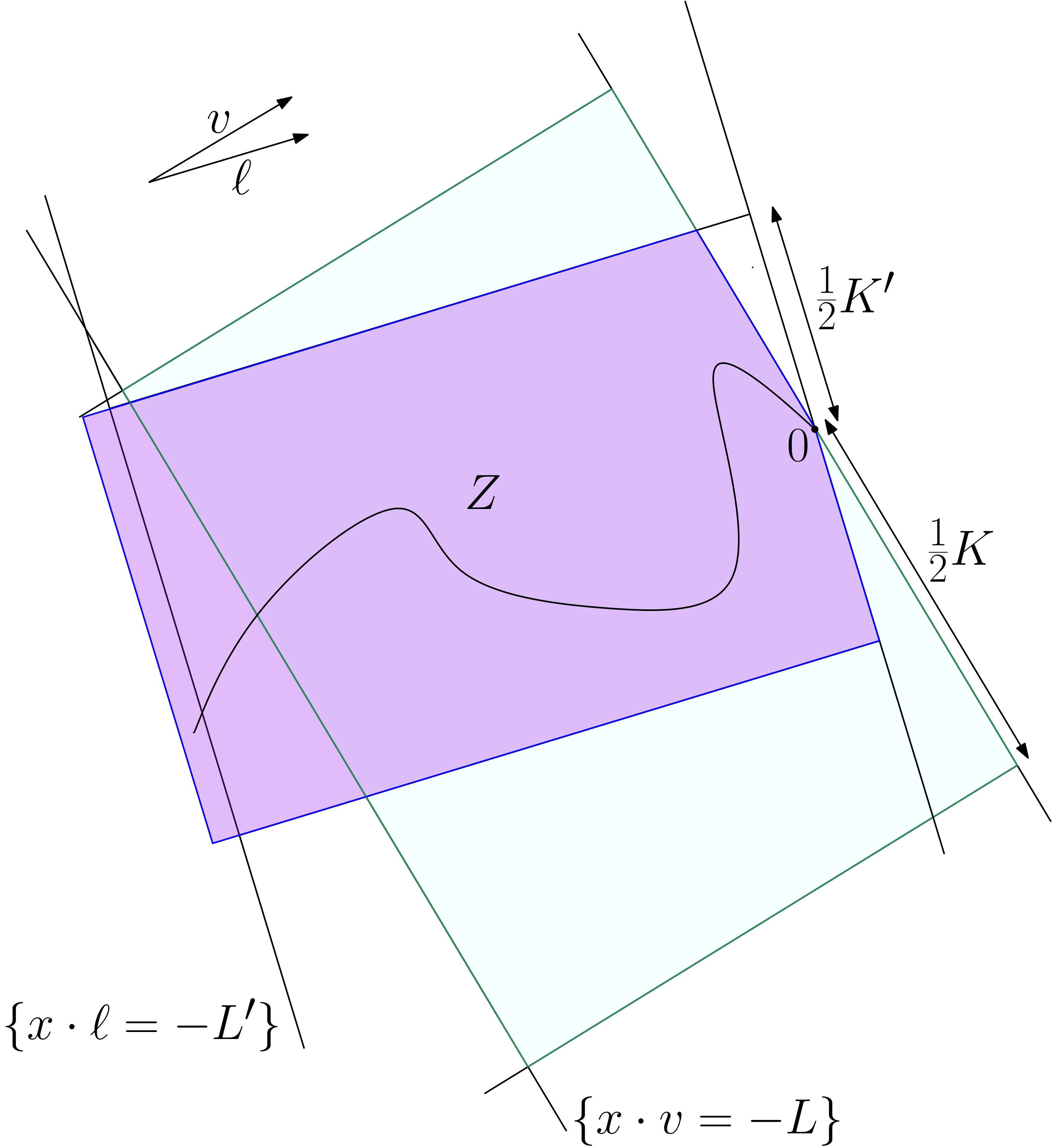}
    \caption{In order for the walk to cross the line $\{x\cdot\ell=-L'\}$ before leaving the set $Z$, it must exit the lighter shaded box through the line $\{x\cdot v=-L\}$.}
    \label{fig:RationalIrrational}
\end{figure}
Now on the event $\{T_{\leq-L'}^{\ell}<\tilde{T}_{\geq0}^{\ell}\wedge \tilde{T}_{\geq0}^{\vecv}\wedge H_{\geq \frac12K'}^{\ell}\}$, it is necessarily the case that $X_n\in Z$ for $0\leq n\leq T_{\leq-L'}^{\ell}$. Furthermore, if $z:=X_{T_{\leq-L'}^{\ell}}$, then since $z\cdot\ell\leq-L'$ and $z\cdot\ell^{\perp}\leq\frac{K'}{2}$ for all $\ell^{\perp}\perp\ell$, the choice of $\vecv$ implies that $z\cdot\vecv\leq-L$. Thus $T_{\leq-L}^v\leq T_{\leq-L'}^{\ell}$, so $X_n\in Z$ for $0\leq n\leq T_{\leq -L}^{\vecv}$, and therefore $T_{\leq-L}^{\vecv}<\tilde{T}_{\geq0}^{\vecv}\wedge H_{\geq\frac12K}^{\vecv}$. Hence (using \eqref{eqn:591}),
$$
\Pa_{\Gr}^0(T_{\leq-L}^{\vecv}<\tilde{T}_{\geq0}^{\vecv}\wedge H_{\geq\frac12K}^{\vecv})\geq\Pa_{\Gr}^0(T_{\leq-L'}^{\ell}<\tilde{T}_{\geq0}^{\ell}\wedge \tilde{T}_{\geq0}^{\vecv}\wedge H_{\geq \frac12K'}^{\ell})>\alpha.
$$
This is \eqref{eqn:594}. 

For $L\geq0$, let $\vecv=\vecv(L)$ and $K=K(L)$ be defined as in \eqref{eqn:594}, with $K$ increasing in $L$. As before, let $\vecu$ be a constant multiple of $\vecv$ such that $\vecu\in\Z^d$and let $(\vecuk,\vecuk_2,\ldots,\vecuk_d)$ be an orthogonal basis for $\R^d$ such that $\vecuk_i\in\Z^d$ for all $i$, and define $N$ as before as well. 

We define the graph $\Hg_{N,L}$ as described before, using the rational direction $\vecv$, rather than the direction $\ell$, to define it. Thus,
$$
C_{N,L}:=\{x\in\Z^d:0\leq x\cdot\vecv\leq L\}/(N\Z\vecu_2,\ldots,N\Z\vecu_d),
$$
We note here our reason for using $L$ as the length of the cylinder, rather than $L|\vecu|$ as in \cite{Tournier2015}. The choice of $\vecv$ depends on the length of the cylinder, but $|\vecu|$ depends on $\vecv$, and may be unbounded as $\vecv\to\ell$.

As before, arguments based on the graph $\Hg_{N,L}$ give us
\begin{equation}\label{eqn:801}
\frac{1}{2}
\leq 
\frac{1}{2}\Pa_{\Hg_{N,L}}^{M}(\tilde{T}_{\partial}= \tau)
+
\Pa_{\Hg_{N,L}}^{\partial}(X_1\neq  M,T_{M}<\tilde{T}_{\partial}),
\end{equation}
and we must show that the term $\Pa_{\Hg_{N,L}}^{\partial}(X_1\neq  M,T_{M}<\tilde{T}_{\partial})$ approaches 0 as $L$ and $K$ increase. Defining $B=B(L,K)$ as before, our previous arguments give us
\begin{align}
\label{eqn:810}
    \Pa_{\Hg_{N,L}}^{\partial}(X_1\neq  M,T_{M}<\tilde{T}_{\partial})
    &=
    \Pa_{\Gr}^0(T_{B^c}<T_{\leq-R}^{\vecvk})
\end{align}

Comparing with \eqref{eqn:669}, the only difference is that the right hand side considers the event $\{T_{B^c}<T_{\leq-R}^{\vecvk}\}$, rather than $\{T_{B^c}<T_{\leq-R}^{\ell}\}$.

We now must show that \eqref{eqn:810} goes to 0 as $L$ increases (along with $N$, and with $\vecuk$ approaching $\ell$). Let $\varepsilon>0$ and choose $R'>R$. Just as $\Pa_{\Gr}^0(T_{\leq -R}^{\ell}<\infty)=1$, we have $\Pa_{\Gr}^0(T_{\leq -R'}^{\ell}<\infty)=1$. Choose $Q=Q(\varepsilon)$ so that $\Pa_{\Gr}^0(T_{Q^c}\leq T_{\leq-R'}^{\ell})<\varepsilon$.
For large enough $L$, as in \eqref{eqn:874}, we have
\begin{equation}\label{eqn:1008}
    \{T_{B^c}\leq T_{\leq-R}^{\vecv}\}\subset\{T_{Q^c}\leq T_{\leq-R}^{\vecv}\}.
\end{equation}
Now by Lemma \ref{lem:598}, for $\vecv$ close enough to $\ell$ (i.e., for large enough $L$), if $x\cdot\ell\leq -R'$ and $x\cdot\vecv\geq -R$, then $x$ is not in $Q$, so that
the event $\{T_{\leq-R'}^{\ell}\leq T_{Q^c}\leq T_{\leq -R}^{\vecv}\}$ is impossible, and therefore 
\begin{equation}\label{eqn:1013}
    \{T_{Q^c}\leq T_{\leq-R}^{\vecv}\}\subset\{T_{Q^c}\leq T_{\leq-R'}^{\ell}\}.
\end{equation}
It follows from \eqref{eqn:1008}, \eqref{eqn:1013}, and the choice of $Q$ that for large enough $L$,
$$
\Pa_{\Gr}^0(T_{B^c}\leq T_{\leq-R}^{\vecv})\leq \Pa_{\Gr}^0(T_{Q^c}\leq T_{\leq-R'}^{\ell})<\varepsilon.
$$
Since this can be true for arbitrary epsilon, $\Pa_{\Gr}^0(T_{B^c}\leq T_{\leq-R}^{\vecv})$ goes to 0, and therefore so does $\Pa_{\Hg_{N,L}}^{\partial}(X_1\neq  M,T_{M}<\tilde{T}_{\partial})$.

Next, we will must show that $\Pa_{\Hg_{N,L}}^M(T_{\partial}=\tau)$ is bounded away from 1 as $M$ increases. Using \eqref{eqn:594} in place of \eqref{eqn:585}, we are able to argue as before to get
\begin{align*}
    \Pa_{\Hg_{N,L}}^M(T_{\partial}\neq \tau)
    &\geq \frac12\alpha.
\end{align*} 

Now taking the limsup in \eqref{eqn:801} as $M\to\infty$ yields the contradiction

$$
\frac{1}{2}
\leq 
\frac{1}{2}\left(1-\frac12\alpha\right)
<
\frac{1}{2}.
$$
\end{proof}

We now have enough to prove Theorem \ref{thm:DirichletFullStatement}.

\begin{proof}[Proof of Theorem \ref{thm:DirichletFullStatement}]\hypertarget{proof:DirichletFullStatement}
First, suppose $\vartriangle\neq0$. Then if $\ell\cdot\vartriangle>0$, the arguments in \cite{Tournier2015}, which work for bounded jumps, show that $\Pa^0(A_{\ell})>0$. For $d\geq3$, the proof of the 0-1 law in \cite{Bouchet2013} can easily be modified to work for bounded jumps, as a remark in \cite{Tournier2015} points out. If $d=1$, the 0-1 law of \cite{Key1984} applies, and if $d=2$, Theorem \ref{thm:01Law} applies. Thus, we get $\Pa^0(A_{\ell})=1$. If $\ell\cdot\vartriangle<0$, then $-\ell\cdot\vartriangle>0$, so we get $\Pa^0(A_{-\ell})=1$, and therefore $\Pa^0(A_{\ell})=0$. Finally, if $\ell\cdot\vartriangle=0$, then the results of \cite{Drewitz&Ramirez2010} (which can easily be made to work for bounded jumps, as noted in the aforementioned remark in \cite{Tournier2015}) imply that $\Pa^0(A_{\ell})=0$. This handles the case $\vartriangle\neq0$.
On the other hand, if $\vartriangle=0$, then the conclusion is that of Theorem \ref{thm:Dirichlet}.
\end{proof}

\subsection{Further remarks}
We have generalized to RWDE with bounded jumps the complete characterization of $\Pa_{\Gr}^0(A_{\ell})$ that was known for nearest-neighbor RWDE.

However, there is one nagging difficulty in the zero-drift case that must be dealt with before we may claim absolute victory over the issue of directional transience for RWDE. Because there are uncountably many directions, proving that the probability of transience in any given direction is zero does not automatically mean that it is impossible for the walk to be directionally transient. For example, however unlikely it seems, one could imagine the possibility that a walk is almost surely transient in some random direction $\ell\in S^{d-1}$ with continuous (or otherwise atom-free) distribution and recurrent in all directions $\ell'\neq\pm\ell$. This pathological behavior has yet to be ruled out, even for the nearest-neighbor Dirichlet case. To resolve this difficulty we would need to prove, at least for Dirichlet environments, a strengthened version of Conjecture \ref{conj1} (recall that Conjecture \ref{conj1} states that for all $\ell\in S^{d-1}$, $\Pa^0A_{\ell}^0)=0$), where $A_{\ell}^0$ is the event that the walk is transient in direction $\ell$, but not in a neighborhood of directions around $\ell$). 
\begin{conj}\label{conj2}
Let $\Pa^0$ be the law of an i.i.d. RWRE on $\Z^d$. Then 
$\Pa^0\left(\bigcup_{\ell\in S^{d-1}}A_{\ell}^0\right)=0$.
\end{conj}

\section*{Acknowledgements}

The author thanks his advisor, Jonathon Peterson, for his mentorship while working on this problem. Part of the work on this paper was supported by a Bilsland dissertation fellowship through Purdue University. The author wrote Appendix \ref{append:C2} while supported by NSF grant DMS-2153869.

\bibliographystyle{plain}
\bibliography{default}

\appendices

\section{Kalikow's 0-1 Law and Other Lemmas}\label{app:Kalikow}

In the proof of Theorem \ref{thm:01Law}, we appealed to the following theorem.
\begin{thm}[Kalikow's 0-1 Law]\label{thm:K01Law}
    Let $\Pa^0$ be the annealed measure of a RWRE on $\Z^d$ satisfying assumptions (\hyperlink{cond:C1}{C1}), (\hyperlink{cond:C2}{C2}), and (\hyperlink{cond:C3}{C3}). Then for every $\ell\in S^{d-1}$, $\Pa^0(A_{\ell}\cup A_{-\ell})\in\{0,1\}$. 
\end{thm}
We also appealed to the following lemma, which is an ingredient in the proof of Theorem \ref{thm:K01Law}.
\begin{lem}\label{lem:986}
$\Pa^0(A_{\ell})>0$ if and only if $\Pa^0(T_{<0}=\infty)>0$
and $\Pa^0(A_{-\ell})>0$ if and only if $\Pa^0(T_{>0}=\infty)>0$.
\end{lem}

A rudimentary version of Theorem \ref{thm:K01Law} for two dimensions, using a uniform ellipticity assumption and assuming $\ell=(0,1)$, was first given by Kalikow in \cite{Kalikow1981}. Improvements were made in \cite{Sznitman&Zerner1999} (allowing general $d$ and general $\ell$) and \cite{Zerner&Merkl2001} (removing the uniform ellipticity assumption), but the overall structure of the argument has changed very little. The proof in \cite{Zerner&Merkl2001} does not use the nearest-neighbor assumption, except in a version of Lemma \ref{lem:NoGettingStuckInaStrip}, which we now state and re-prove using the same ideas but without the nearest-neighbor assumption.

\begin{proof}[Proof of Lemma \ref{lem:NoGettingStuckInaStrip}]
Observe that on the event in question, it is either the case that for some $y$ with $a\leq y\cdot\ell<b$, $X_n=y$ infinitely often, or that $X_n$ hits infinitely many vertices in the slab $\{a\leq x\cdot\ell<b\}$. It therefore suffices to show that the intersection of each of these events with the event $T_{\geq b}^{\ell}=\infty$ has probability 0. 

First, fix $y$  with $a\leq y\cdot\ell<b$. By the irreducibility assumption (\hyperlink{cond:C2}{C2}), $P_{\omega}^y(T_{\geq b}^{\ell}<\tilde{T}_y)>0$ for almost every $\omega$. For such an $\omega$, the strong Markov property implies that the quenched probability of hitting $y$ at least $n$ times before $T_{\geq b}^{\ell}$ is no more than $P_{\omega}^y(\tilde{T}_y<T_{\geq b}^{\ell})^{n-1}$, which approaches 0 as $n\to\infty$. Thus, the (quenched or annealed) probability of hitting $y$ infinitely many times without ever reaching the half-space $\{x\cdot\ell\geq b\}$ is 0. Summing over countably many $y$ still gives a probability of 0. 

Now consider the event that infinitely many points in $\{a \leq x\cdot\ell<b\}$ are hit. By assumption (\hyperlink{cond:C2}{C2}), each of these points $x$ has a {\em possible path} (in the notation introduced in the proof of Theorem \ref{thm:01Law}) to $\{x\cdot\ell\geq b\}$. By shift-invariance, there is some $N>0$ and $\varepsilon>0$ such that each $x$ in $\{a \leq x\cdot\ell<b\}$ has a possible path of length no more than $N$ and with annealed probability at least $\varepsilon$. Thus, in order to hit infinitely many points in $\{a\leq x\cdot\ell< b\}$, the walk must hit the vertex sets of infinitely many disjoint paths to $\{x\cdot\ell\geq b\}$, each of which has length no more than $N$ and annealed probability at least $\varepsilon$. Now by the i.i.d. assumption (\hyperlink{cond:C1}{C1}), each time the walk hits an unexplored vertex set of such a path, its probability, conditioned on its entire past, of immediately taking (the rest of) that path is at least $\varepsilon$. The probability of hitting vertex sets of $n$ unique such paths before hitting $\{x\cdot\ell\geq b\}$ is therefore no more than $(1-\varepsilon)^{n-1}$, which approaches 0 as $n\to\infty$. Thus, the annealed probability of hitting the vertex sets of infinitely many disjoint paths to $\{x\cdot\ell\geq b\}$, each with length no more than $N$ and annealed probability at least $\varepsilon$, without ever reaching $\{x\cdot\ell\geq b\}$, is 0. Since these paths have no more than $N$ vertices in them, hitting infinitely many sites in $\{a\leq x\cdot\ell <b\}$ requires hitting the disjoint vertex sets of infinitely many such paths, and therefore the annealed probability of hitting infinitely many sites in $\{a\leq x\cdot\ell <b\}$ without ever reaching $\{x\cdot\ell\geq b\}$ is 0. This gives us \eqref{eqn:11}.
\end{proof}

Lemma \ref{lem:986} is an easy consequence of Lemma \ref{lem:NoGettingStuckInaStrip} (see \cite{Sznitman&Zerner1999}, \cite{Zerner2007}). We repeat the proof here because it is short.
\begin{proof}[Proof of Lemma \ref{lem:986}]
    By Lemma \ref{lem:NoGettingStuckInaStrip}, the event $\{T_{<0}=\infty\}\setminus A_{\ell}$ must have zero probability, so $\Pa^0(T_{<0}=\infty)>0$ implies $\Pa^0(A_{\ell})>0$. 
    On the other hand, suppose $\Pa^0(T_{<0}=\infty)=0$. By shift-invariance, $\Pa^x(T_{<x\cdot \ell}=\infty)=0$ for all $x\in\Z^2$, which implies that with $\Pm$-probability 1, $P_{\omega}^x(T_{<x\cdot \ell}<\infty)=1$ for all $x$. By the strong Markov property, this implies $\Pa^0(A_{\ell})=0$. 
    For the second result, just take $\ell'=-\ell$.
\end{proof}

Our last lemma justifies the assertion in \eqref{eqn:558}.

\begin{lem}\label{lem:Conditioning}
For all measurable events $A,B\subset(\Z^2)^{\N}$ with $\Pa^a(A)>0$    
$$
\Pa^a(B|A)=\Ea^a[P_{\omega}^a(B|A)|A]
$$
\end{lem}

\begin{proof}
    For any probability measure $P$, random variable $Y$, and event $A$ of positive probability, we have $E[Y|A]=\frac{E[Y\mathbbm{1}_A]}{P(A)}$. Using $\Pa^a$ as our probability measure and letting $Y=P_{\omega}^a(B|A)$, we get
    \begin{align*}
        \Ea^a[P_{\omega}^a(B|A)|A]&=\frac{\Ea^a[P_{\omega}^a(B|A)\mathbbm{1}_A]}{\Pa^a(A)}
        \\
        &=\frac1{\Pa^a(A)}\Ea^a\left[\frac{P_{\omega}^a(B\cap A)}{P_{\omega}^a(A)}\mathbbm{1}_A\right]
        \\
        &=\frac1{\Pa^a(A)}\Ea^a\left[\Ea^a\left[\frac{P_{\omega}^a(B\cap A)}{P_{\omega}^a(A)}\mathbbm{1}_A\given \omega\right]\right].
        \end{align*}
        Now pulling the part that is $\omega$-measurable out of the conditional expectation and noting that $\Ea^a[\mathbbm{1}_A|\omega]$ is simply $P_{\omega}^a(A)$, we have
        \begin{align*}
        \Ea^a[P_{\omega}^a(B|A)|A]
        &=\frac1{\Pa^a(A)}\Ea^a\left[\frac{P_{\omega}^a(B\cap A)}{P_{\omega}^a(A)}\Ea^a\left[\mathbbm{1}_A\given \omega\right]\right]
        \\
        &=\frac1{\Pa^a(A)}\Ea^a\left[\frac{P_{\omega}^a(B\cap A)}{P_{\omega}^a(A)}P_{\omega}^a(A)\right]
        \\
        &=\frac{\Ea^a[P_{\omega}^a(B\cap A)]}{\Pa^a(A)}
        \\
        &=\Pa^a(B|A).
    \end{align*}
\end{proof}

\section{A non-elliptic model satisfying the conditions for Theorem \ref{thm:01Law}
}\label{append:C2}

In this appendix, we present an example of a two-dimensional RWRE model where the Markov chain is almost surely not irreducible, but nevertheless conditions (\hyperlink{cond:C1}{C1}), (\hyperlink{cond:C2}{C2}), and (\hyperlink{cond:C3}{C3}) are satisfied. Recall that (\hyperlink{cond:C1}{C1}) is the i.i.d. condition, condition (\hyperlink{cond:C3}{C3}) requires bounded jumps, and (\hyperlink{cond:C2}{C2}), which replaces the ellipticity condition of \cite{Zerner&Merkl2001}, states that with $\Pm$-probability 1, the Markov chain induced by $\omega$ has only one infinite communicating class, and it is reachable from every site.

\begin{figure}
    \centering
    \includegraphics[width=3in]{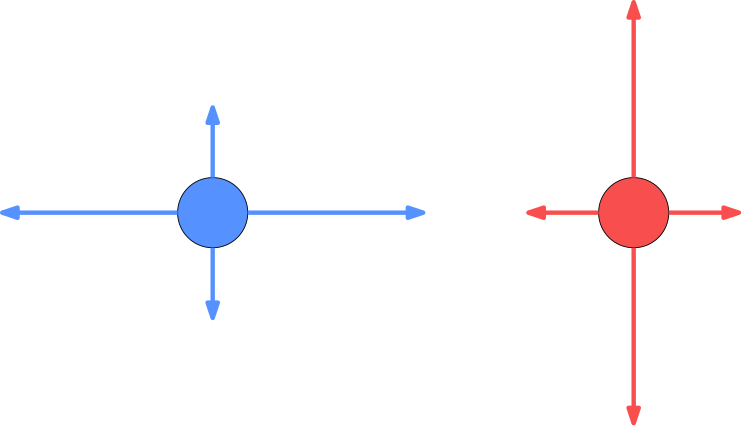}
    \caption{A blue site is pictured on the left and a red site on the right, with their possible jumps.}
    \label{fig:node}
\end{figure}

For an environment $\omega\in\Omega_{\Z^2}$, say a site $x\in\Z^2$ is \textit{blue} under $\omega$ if $\omega(x,x+y)>0\Leftrightarrow y\in\{\pm2e_1,\pm e_2\}$, where $e_1$ and $e_2$ are the sites directly to the right of and above the origin, respectively. Say $x$ is \textit{red} if  $\omega(x,x+y)>0\Leftrightarrow y\in\{\pm e_1,\pm 2e_2\}$. See figure \ref{fig:node}.
Consider an i.i.d. measure $\Pm$ on $\Omega_{\Z^2}$ (that is, a measure satisfying condition (\hyperlink{cond:C1}{C1}) whose marginals are such that each site $x$ is blue with probability $p\in(0,1)$ and red with probability $(1-p)$ (note that conditioned on a site being red or blue, the values of the transition probabilities that are specified to be positive may still be random).

\begin{figure}
    \centering
    \includegraphics[width=2.5in]{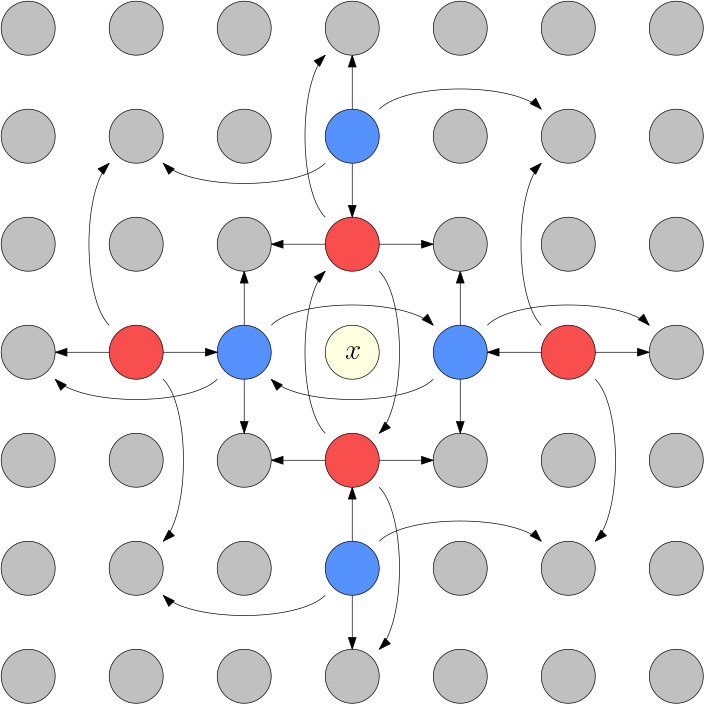}
    \caption{A site $x$ that is not reachable from any other site.}
    \label{fig:NotIrreducible}
\end{figure}

Under $\Pm$, the environment will almost surely not be irreducible. Indeed, if the two nearest vertical neighbors of $x$ are red, and the two sites above and below these are blue, and the nearest horizontal neighbors are blue, and the two sites to the right and left of these are red, then $x$ will not be reachable from any other vertex. See figure \ref{fig:NotIrreducible}. Nevertheless, $\Pm$ satisfies the assumptions needed for Theorem \ref{thm:01Law}.

\begin{prop}\label{prop:FunModel}
    The measure $\Pm$ described in this appendix satisfies conditions (\hyperlink{cond:C1}{C1}), (\hyperlink{cond:C2}{C2}), and (\hyperlink{cond:C3}{C3}).
\end{prop}

We break the proof up into a series of lemmas.

\begin{lem}
    For almost every $\omega$, there exists  $j\geq1$ such that $je_1$ is reachable from every site in $\Z^2$.
\end{lem}

\begin{proof}
    Let $j\geq1$ be the lowest positive integer such that $je_1$ and $(j+1)e_1$ are both red. Such $j$ exists for almost every environment $\omega$. Now from every site, it is possible to step to the right, to the left, up, or down (in this argument, ``possible'' means possible with positive probability in almost every environment $\omega$). Therefore, from every site, it is possible to proceed to a point that is above and to the right of the origin. From there, it is possible to step down repeatedly, either to a site that is horizontally level with the origin (i.e., $ke_1$ for some $k\in\N$), or a site that is one unit higher than the origin (i.e., $e_2+ke_1$) for some $k\in\N$). In the latter case, if one is at a blue site, it is possible to step down to $ke_1$. Otherwise, it is possible to step repeatedly to the right until reaching a blue site, which happens in finitely many steps for almost every $\omega$, and then step down from there to $k'e_1$. Now by stepping repeatedly to the right or to the left, one can reach either $je_1$ or $(j+1)e_1$, and since $(j+1)e_1$ is red, one may step from there to $je_1$. 
\end{proof}

\begin{lem}\label{lem:1212}
    There almost surely exists an infinite communicating class $\mathcal{C}$ which is reachable from every site. 
\end{lem}

\begin{proof}
    Take the set containing $je_1$ and all sites reachable from there. We note that $\mathcal{C}$ is precisely the set of points that are reachable from every site.
\end{proof}

It remains to show that there is no infinite communicating class but $\mathcal{C}$. 

\begin{lem}\label{lem:1222}
    For almost every $\omega$, if a site $x\in\Z^2$ is not in $\mathcal{C}$, then each of its four nearest neighbors is in $\mathcal{C}$. Moreover, the sites $x\pm e_1$ are blue, and the sites $x\pm e_2$ are red.
\end{lem}

\begin{proof}
    We already know $je_1$ is reachable from every site. From there, it is possible to step to the left until hitting the origin or $e_1$. Therefore, if the origin is not in $\mathcal{C}$, then $e_1$ is reachable from every site and thus in $\mathcal{C}$. By symmetry, if the origin is not in $\mathcal{C}$ then every nearest neighbor of the origin is in $\mathcal{C}$. By shift-invariance, this is true of every site $x$.
    Now if a site immediately to the right of $x$ were red or a site immediately to the left were blue, it would be possible to step from such a site to $x$, putting $x$ in $\mathcal{C}$. Therefore, the sites $x\pm e_1$ are blue, and the sites $x\pm e_2$ are red.
\end{proof}

\begin{lem}\label{lem:1232}
    Let $\mathcal{D}=\mathcal{D}(\omega)$ be the communicating class containing the origin. It is almost surely the case that if $\mathcal{D}\neq\mathcal{C}$, then $|\mathcal{D}|<\infty$.
\end{lem}

\begin{figure}
    \centering
    \includegraphics[width=4.5in]{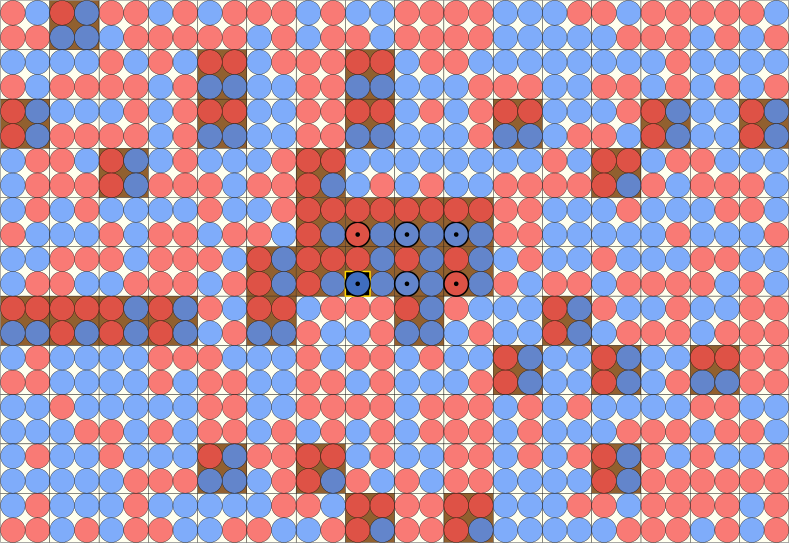}
    \caption{A portion of an environment $\omega$. Squares of four elements, each with the bottom left element in $(2\Z)^2$, are shown with a brown background if the bottom left element is in $\mathcal{B}$. The origin his highlighted with a black background, and vertices in $\mathcal{D}$ are marked with dots.}
    \label{fig:UniqueExtensionSignedMeasure}
\end{figure}

\begin{proof}
    Suppose $\mathcal{D}\notin\mathcal{C}$. Then by Lemma \ref{lem:1222}, any site in $\mathcal{D}$ that is reachable in one step from the origin must almost surely be in the ``even sublattice'' $(2\Z)^2$. By induction, applying Lemma \ref{lem:1222} repeatedly, this conclusion extends to every site in $\mathcal{D}$ that is reachable in any number of steps from the origin. But by definition, every site in $\mathcal{D}$ is reachable from the origin, so $\mathcal{D}\subseteq(2\Z^2)$. Now in order for a site $x\in(2\Z)^2$ to be in $\mathcal{D}$, it must avoid being in $\mathcal{C}$, which by  Lemma \ref{lem:1222} almost surely requires at least that the site to its immediate right be blue and the site immediately above it be red. Let $\mathcal{B}$ be the set of sites $x\in(2\Z)^2$ such that $x+e_1$ is blue and $x+e_2$ is red. Then almost surely, if $\mathcal{D}\neq\mathcal{C}$, then $\mathcal{D}$ is a subset of $\mathcal{B}$ which is connected in the superlattice $(2\Z)^2$.

    Notice, however, that the events $\{x\in\mathcal{B}\}_{x\in2\Z}$ are independent, and each have probability $p(1-p)\leq\frac14$ (recall that $p$ is the probability that a given site is blue). Now $\frac14$ is well below the critical percolation threshhold\footnote{Numerical estimates put this threshold around 0.59274621; see \cite{Newman&Ziff2000}.} for site percolation on $\Z^2$. Indeed, a union bound puts the threshhold at least at $\frac13$, since there are at most $3^n$ loop-free paths of length $n$ from the origin. If the occupation probability is less than $\frac13$, this means the probability that there exists an occupied  path of to a point of distance $n$ from the origin (under, say, the lattice metric) decays exponentially in $n$. Therefore, every component of $\mathcal{B}$ that is ``connected in $(2\Z)^2$'' (meaning any two two points in a component are joined by a nearest-neighbor path in $(2\Z)^2$) is finite. In particular, $\mathcal{D}$ is finite, since it is a subset of such a component.     
\end{proof}

We now have enough to prove that our model satisfies (\hyperlink{cond:C1}{C1}), (\hyperlink{cond:C2}{C2}), and (\hyperlink{cond:C3}{C3}).

\begin{proof}[Proof of Proposition \ref{prop:FunModel}]
    It satisfies condition (\hyperlink{cond:C1}{C1}) by assumption. Likewise, it satisfies condition (\hyperlink{cond:C3}{C3}) with $R=2$. Thus, we need only show that condition (\hyperlink{cond:C2}{C2}) is satisfied. The existence of a communicating class $\mathcal{C}$ that is infinite and reachable from every site is the content of Lemma \ref{lem:1212}. To show that there is only one infinite communicating class, it suffices by shift-invariance to show that if the origin is not in $\mathcal{C}$ then it is in a finite communicating class. This is the content of Lemma \ref{lem:1232}.
\end{proof}

\end{document}